\documentclass[11pt]{article}

\usepackage{amssymb,latexsym,amsmath}

\begin{document}

\newcommand{\bfi}{\bfseries\itshape}

\makeatletter

\@addtoreset{figure}{section}

\def\thefigure{\thesection.\@arabic\c@figure}

\def\fps@figure{h, t}

\@addtoreset{table}{bsection}

\def\thetable{\thesection.\@arabic\c@table}

\def\fps@table{h, t}

\@addtoreset{equation}{section}

\def\theequation{\thesubsection.\arabic{equation}}

\makeatother

\newtheorem{thm}{Theorem}[section]

\newtheorem{prop}[thm]{Proposition}

\newtheorem{lema}[thm]{Lemma}

\newtheorem{cor}[thm]{Corollary}

\newtheorem{defi}[thm]{Definition}

\newtheorem{rk}[thm]{Remark}

\newtheorem{exempl}{Example}[section]

\newenvironment{exemplu}{\begin{exempl}  \em}{\hfill $\surd$

\end{exempl}}

\newcommand{\comment}[1]{\par\noindent{\raggedright\texttt{#1}

\par\marginpar{\textsc{Comment}}}}

\newcommand{\todo}[1]{\vspace{5 mm}\par \noindent \marginpar{\textsc{ToDo}}\framebox{\begin{minipage}[c]{0.95 \textwidth}

\tt #1 \end{minipage}}\vspace{5 mm}\par}

\newcommand{\ea}{\mbox{{\bf a}}}
\newcommand{\eu}{\mbox{{\bf u}}}
\newcommand{\ueu}{\underline{\eu}}
\newcommand{\ueo}{\overline{u}}
\newcommand{\oeu}{\overline{\eu}}
\newcommand{\ew}{\mbox{{\bf w}}}
\newcommand{\ef}{\mbox{{\bf f}}}
\newcommand{\eF}{\mbox{{\bf F}}}
\newcommand{\eC}{\mbox{{\bf C}}}
\newcommand{\en}{\mbox{{\bf n}}}
\newcommand{\eT}{\mbox{{\bf T}}}
\newcommand{\eL}{\mbox{{\bf L}}}
\newcommand{\eV}{\mbox{{\bf V}}}
\newcommand{\eU}{\mbox{{\bf U}}}
\newcommand{\ev}{\mbox{{\bf v}}}
\newcommand{\eve}{\mbox{{\bf e}}}
\newcommand{\uev}{\underline{\ev}}
\newcommand{\eY}{\mbox{{\bf Y}}}
\newcommand{\eK}{\mbox{{\bf K}}}
\newcommand{\eP}{\mbox{{\bf P}}}
\newcommand{\eS}{\mbox{{\bf S}}}
\newcommand{\eJ}{\mbox{{\bf J}}}
\newcommand{\eB}{\mbox{{\bf B}}}
\newcommand{\leb}{\mathcal{ L}^{n}}
\newcommand{\eI}{\mathcal{ I}}
\newcommand{\eE}{\mathcal{ E}}
\newcommand{\hen}{\mathcal{ H}^{n-1}}
\newcommand{\eBV}{\mbox{{\bf BV}}}
\newcommand{\eA}{\mbox{{\bf A}}}
\newcommand{\eSBV}{\mbox{{\bf SBV}}}
\newcommand{\eBD}{\mbox{{\bf BD}}}
\newcommand{\eSBD}{\mbox{{\bf SBD}}}
\newcommand{\ecs}{\mbox{{\bf X}}}
\newcommand{\eg}{\mbox{{\bf g}}}
\newcommand{\paromega}{\partial \Omega}
\newcommand{\gau}{\Gamma_{u}}
\newcommand{\gaf}{\Gamma_{f}}
\newcommand{\sig}{{\bf \sigma}}
\newcommand{\gac}{\Gamma_{\mbox{{\bf c}}}}
\newcommand{\deu}{\dot{\eu}}
\newcommand{\dueu}{\underline{\deu}}
\newcommand{\dev}{\dot{\ev}}
\newcommand{\duev}{\underline{\dev}}
\newcommand{\weak}{\rightharpoonup}
\newcommand{\weakdown}{\rightharpoondown}
\renewcommand{\contentsname}{ }

\renewcommand{\contentsname}{ }

\title{Energy minimizing brittle crack propagation II}

\author{Marius Buliga \\
Institute of Mathematics, Romanian Academy \\
P.O. BOX 1-764, RO 70700\\
Bucure\c sti, Romania\\
{\footnotesize Marius.Buliga@imar.ro}}

\date{15.11.1997}

\maketitle

\tableofcontents

\newpage

\section*{Note to the reader}

This paper is an alternative version of the published paper [Bu3], containing more results obtained in 
[Bu1] and no numerical results. The paper is dated by the last modification of the file. 

\section{Introduction}

\indent

This paper is devoted to the study of quasi-static brittle crack evolution. 
We work under the following assumptions: a linear elastic body, with or without initial 
cracks inside, evolves in a quasi-static manner under an imposed path of boundary displacements. 
During its evolution cracks with unprescribed geometry may appear and/or grow. 

The difficulty of brittle crack propagation problems consists in the nature 
of the main unknown: the crack itself, at various moments in time. The research  
in this field  concerns mainly the constitutive behaviour of a brittle 
material, like the basic paper of Griffith [G]. The 
essential was stated in papers like:  Eshelby [Es], Irwin [I], Gurtin [Gu1], [Gu2].  

In almost all the studies the geometry of the crack is   
prescribed. There are few exceptions, as the papers of Ohtsuka [Oht1---3] or Stumpf \& Le [Stle]. 
  The geometry of the crack  can be prescribed in a strong form, like in the case of a plane 
rectangular or elliptic 
crack which is supposed to remain plane rectangular or elliptic during its growth. 
In a weak form, the geometry of the crack can be prescribed by the assumption that  
the configuration of the body is  2 dimensional and the crack is supposed to have 
only an edge, which is a point. In this case the evolution of the crack is conveniently 
reduced to the evolution of a point. Under these assumptions the geometrical nature of the 
main unknown is obscured. 

A new direction of research in brittle fracture mechanics  begins with the article of 
 Mumford \& Shah  [MS] regarding the problem 
of image segmentation. This problem,  which consists in finding the set of edges of a picture 
and constructing 
a smoothed version of that  picture, it turns to be intimately related to the problem of brittle 
crack evolution. 
In the before mentioned article Mumford and Shah propose the
following variational approach to the problem of image segmentation: 
let $g:\Omega \subset\mathbb{R}^{2} \rightarrow [0,1]$ 
be the original picture, given as a distribution of grey levels (1 is white and 0 is black), 
let $u: \Omega \rightarrow R$ 
be the smoothed picture and $K$ be the set of edges. $K$ represents the set where $u$ has jumps, 
i.e. $u \in C^{1}(\Omega \setminus K,R)$. The pair formed by the smoothed picture $u$ and the set 
of edges $K$ minimizes then the functional: 
$$I(u,K) \ = \ \int_{\Omega} \alpha \ \mid \nabla u \mid^{2} \mbox{ d}x \ + \ 
\int_{\Omega} \beta \  \mid u-g \mid^{2} \mbox{ d}x \ + \ \gamma \mathcal{ H}^{1}(K) \ \ .$$
The parameter $\alpha$ controls the smoothness of the new picture $u$, $\beta$ controls the 
$L^{2}$ distance between the smoothed picture and the original one and $\gamma$ controls the 
total length of the edges given by this variational method. The authors remark that for $\beta 
= 0$ the functional $I$ might be useful for an energetic treatment of fracture mechanics. 

An energetic approach to fracture mechanics is naturally suited to explain brittle crack 
appearance under imposed boundary displacements. 
The idea is presented  in the followings. 

The state of a brittle body is  described by a pair 
displacement-crack. $(\eu,K)$ is such a pair if $K$ is a crack --- seen as a surface 
---  
which appears in the body and  $\eu$ is a displacement of the broken body under the imposed 
boundary displacement, i.e. $\eu$ is continuous in the exterior of the surface $K$ and $\eu$
equals the imposed displacement $\eu_{0}$ on the exterior boundary of the body. 

Let us suppose 
that the total energy of the body is a Mumford-Shah functional of the form: 
$$E(\eu,K) \ = \ \int_{\Omega} w(\nabla \eu) \mbox{ d} x \ + \ F(\eu_{0},K) \ \ .$$
The first term of the functional $E$ represents the elastic energy of the body with the
displacement $\eu$. The second term represents the energy consumed to produce the crack 
$K$ in the body, with the boundary displacement $\eu_{0}$   as parameter. 
Then the crack that appears is supposed to be the second term of the pair $(\eu,K)$ which 
minimizes the total energy $E$. 

We shall  use this idea in the study of quasi-static brittle crack evolution. For doing this, 
we proceed to a time discretization which transforms the problem of crack evolution into 
a sequence of energy minimization problems.  Francfort 
\& Marigo [Fma] proceed in the same way in the case of  brittle brutal damage evolution.  However, it 
is only a belief that when the time step goes to zero, the discretized evolution converges to 
an almost continuous (with respect to time) evolution. We have found in the frame of 
generalized minimizing movements, introduced by De Giorgi [DG], stronger mathematical reasons 
to support this belief. That is why we introduce the notion of energy minimizing movement as a 
particular case of a generalized minimizing movement. 

In section 2. the notion of energy minimizing movement is introduced in a form useful in the
sequel. After preliminaries concerning the statics of a brittle body, the Griffith criterion 
of brittle crack propagation is presented in subsection 3.3., as a selection criterion amongst all 
possible crack evolutions. At the end of this section we formulate 
the problem of  quasi-static brittle crack evolution  in the form (\ref{fracpro}). 

In section 4. we give an energy minimizing movement formulation to this problem by 
using a Mumford-Shah energy functional (definition 4.1.). We are assured about the existence of the 
discretized (or incremental) 
solution to that problem by theorem 4.1. (in a weak form presented in the section of proofs). 
In subsection 4.2. we investigate the features of  this first model. In this model 
we have only one material constant connected to fracture, namely the constant of Griffith $G$. 
 We find exact solutions 
and useful estimations in the case of anti-plane displacements (theorem 4.2.), which tell us 
that crack
appearance is allowed in this model. A part of this results can be found in [Bu2], in 
connection with fiber-matrix debonding in composites.  We prove also a bad feature of the model:  
the critical 
stress which lead to fracture in an uni-dimensional  traction experiment is not a constant of 
material. 

In section 5. is presented an improved model, based on a preliminary study of smooth brittle 
crack propagation (see also [Bu1], [Bu3]). The Griffith criterion is reformulated  by using 
proposition 5.1. and the $K2$ functional (definition 5.3.). 
The $K2$ functional is a generalized version of the  J integral of Rice [R].  
After  we modify the Griffith criterion by the extension of  the functional $K2$, we obtain the 
differential global criterion of brittle crack appearance (DA). A stronger version of (DA) is 
the local crack appearance criterion (LA). In subsection  5.4. the improved model is presented. 
In this model we have two  constants of material connected to fracture: $G$ and a quantity with  
dimension 
of stress named $\Sigma$. The critical stress which lead to fracture is deduced from the the 
elastic constants and $\Sigma$, hence this time it is a constant of material.     

Section 6. is devoted to a brief introduction to special functions with bounded variation 
or deformation. Weak versions of theorems 4.1. and 5.2., concerning the existence of the discretized 
(or incremental) solutions of the models presented here,  are given as consequences of 
more general results due to De Giorgi \& Ambrosio [DGA], Ambrosio [A1---3], Belletini, 
Coscia \& Dal Maso, [BCDM].    

In section 7. are given the conclusions regarding the features of the two models presented in the 
paper. We prove a general existence result of the energy minimizing movement described in 
the first model under the assumption of uniformly bounded power communicated by the rest of the 
universe to the body. A  comparation is made with the  model 
of Ambrosio \& Braides [AB], based also on generalized minimizing movements, 
where the crack appearance is forbidden.

\section{General energy minimizing movements}
\indent

An energy minimizing movement is a particular case of a generalized minimizing movement. The 
latter notion has been introduced by De Giorgi in [DG], inspired by the paper [ATW] of Almgren, 
Taylor \& Wang. The definition of a generalized minimizing movement is (according  to Ambrosio
[Amb]) the following: 

\vspace{.5cm} 

{\bf Definition 2.1.} {\it Let $S$ be a topological space and 
$$F: (1,+\infty) \times N \times S \times S \rightarrow\mathbb{R}\cup \left\{ + \infty \right\}$$ 
be a function. For any $u_{0} \in S$, a function $u: [0,+\infty) \rightarrow S$ is a 
 generalized minimizing movement associated to $F$ with 
initial data $u_{0}$,  if  there exists a diverging sequence 
$(s_{i})_{i \in N}$, $s_{i} > 1$, and there are functions $u_{i}: N \rightarrow S$ such that: 

i) $u_{i}(0) \ = \ u_{0}$; 

ii) for any $k \in N$ and any $i$,  $u_{i}(k+1)$ minimizes the functional 
$$v \mapsto F(s_{i}, k ,v, u_{i}(k))$$ over $S$;

iii) for any $t \geq 0$, $u_{i}([s_{i}t]) \rightarrow u(t)$ in $S$ as $i \rightarrow +\infty$ . } 

\vspace{.5cm} 

The canonical example of (generalized) minimizing movement is given by the choice: $S=R^{n}$, 
$f:R^{n} \rightarrow R$ Lipschitz continuous and $C^{2}$ and  
$$F(s,k,u,v) \ = \ f(u) \ + \ \frac{s}{2} \mid u - v \mid^{2} \ \ .$$
In this case, for any $u_{0} \in\mathbb{R}^{n}$ there is only one minimizing movement, 
namely the unique solution of the Cauchy problem 
$$u'(t) \ = \ - \nabla f(u(t)) \ \ \ \  \ , u(0) \ = \ u_{0} \ \ .$$

An energy minimizing movement is a generalized minimizing movement associated to a particular 
function $F$. It is designed to be a "weak stable" solution of  an evolution problem of the
following type: 
\begin{equation}
\left\{ \begin{array}{ll}
\mbox{{\bf A}} \left( u(t),\alpha(t),t \right) \ = \ 0 & \forall \ t \geq 0 \\
\frac{d}{dt} \alpha(t) \ \leq \ \mbox{ {\bf L}} \left( \alpha(t),u(t) \right) &  \forall \ t \geq 0 \\
u(0) \ = \ u_{0} \ \ \ , \alpha(0) \ = \ \alpha_{0} & \ \ \ \ .   
\end{array} \right.
\label{problem}
\end{equation}
 
There are two unknowns in this problem: $u$ and $\alpha$. The evolution of these unknowns is 
quasi-static with respect to $u$. Suppose that we are not in position to give a proper law of
evolution of $\alpha$, or that the law of evolution that we have gives too many solutions. 
Assume further that we have instead at our disposal the expression of the total energy 
of the system described by the pair  $(u,\alpha)$, name it $f(u,\alpha)$, and a set of 
constraints, not in a differential form, upon the evolution of $\alpha$. We can make then a time
discretization with time step $\delta$ and recursively find 
$(u_{k+1}^{\delta},\alpha_{k+1}^{\delta})$ from $(u_{k}^{\delta},\alpha_{k}^{\delta})$, by a
minimization process of the total energy $f$ under some constraints. A weak
stable solution of the previous problem  will be a limit of sequences
$(u_{k}^{\delta},\alpha_{k}^{\delta})_{k}$ 
when the time step $\delta$  converges to $0$. 

In the further definition $S$ may be seen as the space of all 
pairs $x=(u,\alpha)$, endowed with a topology. 

\vspace{.5cm} 

{\bf Definition 2.2.} {\it Let $S$ be a topological space and 

$$F: (1,+\infty) \times N \times S \times S \rightarrow\mathbb{R}\cup \left\{ +\infty \right\} \ \ ,$$
$$F(s,k,x,y) \ = \ f(s,x,y) \ + \ \psi(k/s,y) $$ 
be a function, with $f: N \times S \times S \rightarrow R$ and 
$\psi: [0,\infty) \times S \rightarrow  \left\{ 0, +\infty \right\}$.  For any $x_{0} \in S$, 
a  generalized minimizing movement $x: [0,+\infty) \rightarrow S$ associated to $F$  with 
initial datum $x_{0}$ is an energy minimizing movement associated to the energy $f$ with 
the constraints $\psi$ and initial datum $x_{0}$. }

\vspace{.5cm}

Let us denote by $S(\lambda)$ the following set: 

$$S(\lambda) \ = \ \left\{ y \in S \mbox{ : } \psi(\lambda,y) \ = \ 0 \right\} \ \ .$$ 

From  definition 2.2. we see that $x: [0,\ + \infty) \rightarrow S$ is an energy minimizing
evolution associated to $f$, with the constraints $\psi$ and initial data $x_{0}$ if  there exists 
a diverging sequence $(s_{i})_{i \in N}$ , $s_{i} > 1$, and there are functions 
$x_{i}: N \rightarrow S$  such that: 

i) $x_{i}(0) \ = \ x_{0}$ ;

ii) for any $k \in N$ and any $i \in N$, $x_{i}(k+1)$ minimizes the functional $f$ over the set 
$S(k/s_{i})$ (in particular $x_{i}(k+1)$ belongs to $S(k/s_{i})$); 

iii) for any $t>0$, $x_{i}([s_{i}t]) \rightarrow x(t)$ in $S$ as $i \rightarrow +\infty$ .

\section{Notations and preliminaries}
\indent

\subsection{Notations and constitutive assumptions}
\indent

The open bounded set $\Omega \subset\mathbb{R}^{3}$ represents the reference configuration of an elastic body and 
$\eu:\Omega \rightarrow\mathbb{R}^{3}$ is the displacement field of the body with respect to 
this configuration. We shall always suppose, without mentioning further, 
 that the open set  $\Omega$ and its closure have the same topological boundary. 
 
The expression of the elastic (or free) energy of the body is: 
$$\int_{\Omega} w(\nabla \eu) \mbox{ d}x \ \ .$$ 
The first Piola-Kirchoff stress tensor $\eS$ is 
$$\eS(\eu) \ = \  \frac{d w}{d \nabla} (\nabla \eu)$$ 
and the equilibrium equation of the body in the absence of volumic forces is 
$$div \ \eS(\eu) \ = \ 0 \mbox{ in } \Omega \ \ .$$ 
In this paper we consider that the body is linear elastic and homogeneous, i.e. the 
function $w(\eu)$ has the form: 
$$w(\eu) \ = \ \frac{1}{2}\eC \nabla \eu : \nabla \eu \ \ ,$$
with the elasticity 4-tensor $\eC$ having the symmetries: 
$$\eC_{ijkl} \ = \ \eC_{jikl} \ = \ \eC_{klij} \ \ .$$ 
Under these assumptions the stress tensor $\eS$ becomes the Cauchy stress tensor: 
$$\sig \ = \ \sig(\eu) \ = \ \eC \nabla \eu \ = \ \eC \epsilon(\eu) \ \ ,$$
where $\epsilon(\eu)$ is the symmetric part of $\nabla \eu$, i.e. 
$$\epsilon(\eu) \ = \ \frac{1}{2} \left( \nabla \eu \ + \ \left( \nabla \eu \right)^{T} \right) \ 
\ . $$
 
We shall  suppose moreover that $w$ satisfies the growth conditions: 
$$\forall \eF \in\mathbb{R}^{9} \ , \  \eF \ = \ \eF^{T} \ , \  c \mid \eF \mid^{2} \ \leq \  w(\eF) \ \leq \  C \mid \eF \mid^{2} \ \ ,$$
where $c$ and $C$ belong to $(0,+\infty)$.

In the case of plane displacement  the domain $\Omega \subset\mathbb{R}^{2}$ represents a section in the 
cylindrical reference configuration of the body $\Omega \times R$ and $\eu: \Omega \rightarrow\mathbb{R}^{2}$ 
is a plane displacement. The displacement with respect to the 3 dimensional configuration of the
body has the following  expression: 
$$(x_{1},x_{2},x_{3}) \ \in \ \Omega \times\mathbb{R}\mapsto  (u_{1}(x_{1},x_{2}),u_{2}(x_{1},x_{2}),0) \in\mathbb{R}^{3} \ \ .$$  
In this case we suppose that the body is linear elastic, homogeneous and isotropic. 

In the case of anti-plane displacements the domain $\Omega \subset\mathbb{R}^{2}$ represents a section in the 
cylindrical reference configuration of the body $\Omega \times R$ too. The anti-plane displacement 
is a function $u: \Omega \rightarrow R$. The 3 dimensional displacement has the following form: 
 $$(x_{1},x_{2},x_{3}) \ \in \ \Omega \times\mathbb{R}\mapsto  (0,0,u(x_{1},x_{2})) \in\mathbb{R}^{3} \ \ .$$ 
In this case we make the same assumption of isotropic body, therefore the elastic energy takes the 
form: 
$$\int_{\Omega} \mu \mid \nabla u \mid^{2} \mbox{ d}x \ \ ,$$
where $\mu$ is one of the two Lam\'e's constants.

\subsection{Statics of a fractured elastic body}
\indent

For any measurable set $B \subset\mathbb{R}^{n}$, $\mid B \mid = \leb (B)$ denotes the Lebesgue measure of $B$ and $\mathcal{ H}^{k} (B)$
denotes the $k$ dimensional Hausdorff measure of $B$.

By a crack set in the body $\Omega$ we mean (according with Ball [Ba]) a topologically closed countably 
rectifiable set, generically denoted by $K$. 

Given the function $f$, a point $x \in \Omega$ and an unitary vector (or direction) $\en \in\mathbb{R}^{n}$, the 
approximate limit of $f$ in $x$ with respect to the direction $\en$ is denoted by
$\tilde{f}(x,\en)$ and  it is defined by the following expression: 
$$\lim_{\rho \rightarrow 0_{+}} \frac{\int_{ B_{\rho}(x) \cap \left\{ y \mbox{ : } (y-x)\cdot \en
\geq 0 \right\}} \mid f(y) -  \tilde{f}(x,\en) \mid \mbox{ d}y }
{\mid B_{\rho}(x) \cap \left\{ y \mbox{ : } (y-x)\cdot \en
\geq 0 \right\} \mid} = 0 \ \ .$$
Whenever a field of normals at $K$ is chosen, the lateral 
limits  $f^{+}$ and $f^{-}$ 
of any function $f: \Omega \setminus K \rightarrow\mathbb{R}^{n}$ are $f^{+}:K \rightarrow R$ and  
 $f^{-}:K \rightarrow R$, defined by 
$$f^{+}(x) \ = \ \tilde{f}(x,\en(x)) \ , \ \ f^{-}(x) \ = \ \tilde{f}(x,-\en(x)) \ \ .$$ 
This means that $f^{+}$ and $f^{-}$  satisfy the equalities: 
$$\forall x \in K \ , \ \ 
 \lim_{\rho \rightarrow 0_{+}} \frac{\int_{ B_{\rho}(x) \cap \left\{ y \mbox{ : } (y-x)\cdot \en
\geq 0 \right\}} \mid f(y) -  f^{+}(x) \mid \mbox{ d}y }
{\mid B_{\rho}(x) \cap \left\{ y \mbox{ : } (y-x)\cdot \en
\geq 0 \right\} \mid} = 0 \ \ ,$$
$$\forall x \in K \ , \ \ 
 \lim_{\rho \rightarrow 0_{+}} \frac{\int_{ B_{\rho}(x) \cap \left\{ y \mbox{ : } (y-x)\cdot \en
\leq 0 \right\}} \mid f(y) -  f^{-}(x) \mid \mbox{ d}y }
{\mid B_{\rho}(x) \cap \left\{ y \mbox{ : } (y-x)\cdot \en
\leq 0 \right\} \mid} = 0 \ \ .$$  
Remark that for any $x \in K$ the object $(f^{+}(x), f^{-}(x),\en(x))$ it is unique to a change of 
signs, i.e. 
$$(f^{+}(x), f^{-}(x),\en(x)) \ \sim \ (f^{-}(x), f^{+}(x), -\en(x)) \ \ .$$

For a given crack set $K$ in $\Omega$, by an admissible displacement with respect to $K$ we mean 
a function $\eu: \overline{\Omega} \setminus K \rightarrow\mathbb{R}^{k}$ (where $k$ might be 1,2 or 3) 
which is $C^{1}$  and posses continuous lateral limits on $K$. In this section we shall  consider 
 the space  $W^{1,2}(\Omega \setminus K) \cap 
L^{\infty}(\Omega)$ as the set of weak admissible displacements with respect to the crack set $K$.   

Let $n$ be the dimension of the reference configuration $\Omega$. 
For a given $u_{0} \in H^{\frac{1}{2}}(\paromega,R^{n})\cap L^{\infty}(\paromega,R^{n})$ and for 
a given rectifiable crack set $K$,  such that $\hen(\paromega \setminus K) \ > \ 0$, 
the following problem has a solution $\eu \ = \ \eu(\eu_{0},K)$, unique to rigid displacements 
of $\Omega \setminus K$ equals to 0 on $\paromega$: 
\begin{equation}
\left\{ \begin{array}{ll}
div \ \sig(\eu) \ = \ 0 & \mbox{ in } \Omega \setminus K \\
\sig^{+}(\eu) \en \ = \ \sig^{-}(\eu) \en  \ = \ 0 & \mbox{ on } K \\
\eu \ = \ \eu_{0} & \mbox{ on } \paromega \setminus K \ \ \ .
\end{array}
\right.
\label{echili}
\end{equation}
We use the same notation --- $u=u(u_{0},K)$ --- in the anti-plane case, when $n=2$, $k=1$ 
and the problem (\ref{echili}) becomes 
\begin{equation}
\left\{ \begin{array}{ll}
 \mu \ div \ \nabla u \ = \ 0 & \mbox{ in } \Omega \setminus K \\
(\nabla u)^{+} \en \ = \ (\nabla u)^{-} \en  \ = \ 0 & \mbox{ on } K \\
u \ = \ u_{0} & \mbox{ on } \paromega \setminus K \ \ \ .
\end{array}
\right.
\label{antiechil}
\end{equation}

The solution $\eu(\eu_{0},K)$ of the problem (\ref{echili}) minimizes the functional 
$$E(\ev) \ = \  \int_{\Omega} w(\nabla \ev) \mbox{ d}x$$ 
over the following set of weak  admissible displacements with respect to the crack set $K$: 
$$\left\{ \ev \in W^{1,2}(\Omega \setminus K,R^{n}) \cap L^{\infty}(\Omega,R^{n}) \mbox{ : } 
\ev \ = \ \eu_{0} \mbox{ on } \paromega \setminus K \right\} \ \ .$$
By standard arguments it follows that the functional 
$$\ev \in W^{1,2}(\Omega,R^{n}) \ \mapsto \ \int_{\Omega} \sig(\eu(\eu_{0},K)): 
\nabla \ev \mbox{ d}x$$ 
depends only on the trace of $\ev$ on $\paromega$, hence it give raise to the linear continuous 
function: 
$$\eT(K) : H^{\frac{1}{2}}(\paromega,R^{n}) \cap L^{\infty}(\paromega,R^{n}) \rightarrow  
H^{-\frac{1}{2}}(\paromega,R^{n})  \ \ ,$$ 
\begin{equation}
\langle \eT(K) \eu_{0} , \ev \rangle \ = \  \int_{\Omega} \sig(\eu(\eu_{0},K)): 
\nabla \ev' \mbox{ d}x \ \ \mbox{ for any } \ev' = \ev \mbox{ on } \paromega \ \ .
\label{dirineu}
\end{equation}
In the latter definition $\langle \cdot , \cdot \rangle$ is the duality product of the pair 
of spaces $H^{\frac{1}{2}}(\paromega,R^{n})$ and $H^{-\frac{1}{2}}(\paromega,R^{n})$. 
The function $\eT(K)$ is called the Dirichlet-to-Neumann map of the elastic body $\Omega$ with the 
crack set $K$. 

Under the assumptions concerning the elastic energy density $w$, more precise because of the
symmetries of the elasticity tensor $\eC$, the function $\eT(K)$ is also self-adjoint, i.e. for 
any $\eu, \ev$ we have 
$$\langle \eT(K) \eu,\ev \rangle \ = \ \langle \eT(K) \ev, \eu \rangle \ \ .$$

In the same way  can be defined the Dirichlet-to Neumann map associated to the problem 
(\ref{antiechil}).  

Remark finally that, under the assumptions considered for $w$, the elastic energy of the body 
 can be expressed using the Dirichlet-to-Neumann 
map. Indeed, we have: 
$$\int_{\Omega} w(\nabla \eu(\eu_{0},K)) \mbox{ d}x \ = \ \frac{1}{2} \ \langle \eT(K)
\eu_{0},\eu_{0} \rangle \ \ .$$

\subsection{The Griffith criterion of brittle crack propagation}
\indent

Let us consider in the elastic body $\Omega$ an initial crack set $K_{0}$ which evolves and becomes 
at the moment $t$ the crack set $K_{t}$. We assume that the crack set always increase in time, i.e. 
\begin{equation}
\forall 0 \ < \ t \ < \ t'  \ \ , \ \ K_{t} \ \subset \ K_{t'}  \ \ .
\label{cevol}
\end{equation}

We suppose that the evolution of the body is quasi-static. At the moment $t$ the state of the body 
is characterized by the pair $(\eu(t), K_{t})$, where $\eu(t)$ is the displacement of the body,
admissible with respect to $K_{t}$. Let us denote by $\eu_{0}(t)$ the trace of $\eu(t)$ on 
$\paromega$. We have then the equality $\eu(t) \ = \ \eu(\eu_{0}(t), K_{t})$. 

The power given to the body by the rest of the universe at the moment $t$ has the following
expression: 
$$P(t) \ = \  \int_{\paromega} \eS(\eu(t))\en \cdot \deu_{0}(t) \mbox{ d}x \ = \ 
\langle \eT(K_{t})\eu_{0}(t),\deu_{0}(t) \rangle \ \ .$$ 
Let us consider a given curve $t \mapsto (\eu(t),K_{t})$, such that for any $t$ we have 
$\eu(t) \ = \ \eu(\eu_{0}(t), K_{t})$. For a given $t$ we introduce the following curve 
of displacements: 
$$\forall \tau \geq 0, \ \ \ew(\tau) \ = \ \eu(\eu_{0}(t+\tau), K_{t}) \ \ .$$ 
$\ew(\tau)$ represents the displacement  the body at the moment $t+\tau$ in the presence of 
the crack $K_{t}$. An easy calculation lead us to the equality: 
\begin{equation}
\frac{d}{d\tau}  \int_{\Omega} w(\nabla \ew(\tau)) \mbox{ d}x_{|_{\tau=0}} \ = \ P(t) \ \ . 
\label{extp}
\end{equation}
Therefore $P(t)$ represents the power  consumed at the moment $t$ 
 by the body  in order to modify its displacement, 
constrained to follow the path of imposed boundary displacements  $t \mapsto \eu_{0}(t)$, 
without  any modification of the actual crack set $K_{t}$.    

The Griffith criterion of brittle crack propagation asserts that during the propagation of 
the crack $K_{t}$ the following inequality is true at any moment $t$: 
\begin{equation}
\frac{d}{dt} \left\{ \int_{\Omega} w(\nabla \eu(t)) \mbox{ d}x \ + \ G \hen(K_{t}) \right\} \
\leq \ P(t) \ \ .
\label{grif}
\end{equation}
Here $G$ is the constant of Griffith, supposed to be a material constant. 

The relation (\ref{grif}) can be written in a different form using the Dirichlet-to-Neumann map 
$\eT(K_{t})$. Let us assume that the crack evolution is smooth in the sense that the function 
$t \mapsto \eT(K_{t})$ is differentiable, i.e. the Dirichlet-to-Neumann map varies smoothly in time. 
The Griffith criterion takes the following form: 
$$\frac{1}{2}  \langle \frac{d}{dt} \left[ \eT(K_{t}) \right] \eu_{0}(t),\eu_{0}(t) \rangle \ + \ 
\frac{1}{2} \langle \eT(K_{t}) \deu_{0}(t),\eu_{0}(t) \rangle \ + \ $$ 
$$+ \ \frac{1}{2} \langle \eT(K_{t}) \eu_{0}(t),\deu_{0}(t) \rangle \ + \ G \frac{d}{dt}\left\{
\hen(K_{t}) \right\} \ \leq \ \langle \eT(K_{t})\eu_{0}(t),\deu_{0}(t) \rangle \ \ .$$
The function $\eT(K_{t})$ is self-adjoint, therefore we obtain the following expression of 
the Griffith criterion:
\begin{equation}
 \frac{1}{2}  \langle \frac{d}{dt} \left[ \eT(K_{t}) \right] \eu_{0}(t),\eu_{0}(t) \rangle \ + \
G \frac{d}{dt}\left\{
\hen(K_{t}) \right\} \ \leq \ 0 \ \ .
\label{grifnew}
\end{equation}
We can see that we have the following equality: 
$$P(t) \ - \ \frac{d}{dt}  \int_{\Omega} w(\nabla \eu(t)) \mbox{ d}x \ = \ 
 - \frac{1}{2}  \langle \frac{d}{dt} \left[ \eT(K_{t}) \right] \eu_{0}(t),\eu_{0}(t) \rangle \ \ 
.$$ 
The quantity from the left of the previous equality is usually called the energy release rate 
due only to crack propagation. 

It is  obvious that $\eu_{0}(t)$ plays the role of a time-dependent parameter, since in the last 
inequality $\deu_{0}(t)$ does not appear. 

The problem of quasi-static brittle propagation of an initial crack in an elastic body under a 
time-dependent imposed displacement $\eu_{0}(t)$ can be formally put in the form (\ref{problem}). 
If we put apart the constraint (\ref{cevol}), we have the following formulation: 
\begin{equation}
\left\{ \begin{array}{ll}
\eu(t) \ - \ \eu(\eu_{0}(t),K_{t}) \ = \ 0 & \forall \  t \geq 0 \\
  \frac{1}{2}  \langle \frac{d}{dt} \left[ \eT(K_{t}) \right] \eu_{0}(t),\eu_{0}(t) \rangle \ + \
G \frac{d}{dt}\left\{
\hen(K_{t}) \right\} \ \leq \ 0 &  \forall \  t \geq 0 \\
\eu(0) \ = \ \eu_{0} \ , \ \ K_{0} \ = \ K  \ \ . &  
\end{array} \right.
\label{fracpro}
\end{equation}

\section{The first model}

\indent

In the  left term of the Griffith  criterion (\ref{grif}) appears the time-derivative of an energetic 
functional. Let us consider the following set of admissible pairs displacement-crack: 
$$M \ = \ \left\{ (\eu,K) \mbox{ : } K \mbox{ is a crack set and } \eu \in C^{1}(\overline{\Omega} 
\setminus K,\mathbb{R}^{n}) \mbox{ such that } \right.$$ 
$$\left. (\eu^{+}, \eu^{-}, \en) \mbox{ exists on } K \right\} \ \
.$$ 
The Mumford-Shah energy functional over $M$ has the following expression: 
\begin{equation}
I: M \rightarrow\mathbb{R}\cup \left\{ + \infty \right\} \ , \ \ I(u,K) \ = \ 
\int_{\Omega} w(\nabla \eu) \mbox{ d}x \ + \ G \hen(K) \ \ .
\label{ems}
\end{equation}

\subsection{Introduction of the first model}
\indent

It is natural to try to give an energy minimizing movement formulation to the problem 
(\ref{fracpro}) by using the  functional defined at  (\ref{ems}). According to definition 
2.2. and the constraint (\ref{cevol}), we give the following formulation:   
    
\vspace{.5cm}

{\bf Definition 4.1.} {\it Let us consider the space $M$ endowed with the topology given by the
convergence: 
$$(\eu_{h},K_{h}) \ \rightarrow \ (\eu,K) \mbox{ if } 
\left\{ \begin{array}{l} 
\eu_{h} \ L^{2} \rightarrow \ \eu \\
\hen(K_{h} \Delta K) \ \rightarrow \ 0 \ \ . 
\end{array} \right.$$
We define the functions 
$$J: M \times M \rightarrow\mathbb{R}\ , \ \ J\left( (\eu,K),(\ev,L)\right) \ = \ 
\int_{\Omega} w(\nabla \ev ) \mbox{ d}x \ + \ G \hen(L \setminus K) \ \ ,$$
$$\Psi: [0,\infty) \times M \rightarrow \left\{0,+\infty \right\} \ , \ \ \Psi(\lambda, (\ev,K)) 
\ = \ \left\{ \begin{array}{ll} 
0 & \mbox{ if } \ev = \eu_{0}(\lambda) \mbox{ on } \paromega \setminus K \\
+\infty & \mbox{ otherwise } \ \ . 
\end{array} \right.$$
We consider the initial data $(\eu_{0},K) \in M$ such that $\eu_{0} \ = \ \eu(\eu_{0}(0),K)$. 

For any $s \geq 1$ we recursively define  $(\eu^{s}, K^{s}): N \rightarrow M$ like this: 

\hspace{1.cm} i) $(\eu^{s}, K^{s})(0) \ = \ (\eu_{0},K)$ ;

\hspace{1.cm} ii) for any $k \in N$  $(\eu^{s}, L^{s})(k+1) \in M$  minimizes the functional 
$$ (\ev,L) \in M \ \mapsto \ J\left( ( (\eu^{s},K^{s})(k), (\ev,L) \right) \ + \ \Psi((k+1)/s,
(\ev,L))$$ 
over $M$. In order to verify the constraint (\ref{cevol}), $K^{s}(k+1)$ is defined by the formula: 
$$K^{s}(k+1) \ = \ K^{s}(k) \cup L^{s}(k+1) \ \ .$$

An energy minimizing movement associated to $J$ with the constraints (\ref{cevol}), $\Psi$ and 
initial data $(\eu_{0},K)$  is 
any $(\eu,K): [0,+\infty) \rightarrow M$ having the property:  
there is a diverging sequence $(s_{i})$ such that for any $t > 0$ 
$$(\eu^{s_{i}},K^{s_{i}})([s_{i}t]) \ \rightarrow (\eu,K)(t) \  \mbox{ as } i  \rightarrow \infty \
\ .$$}

\vspace{.5cm}

In the previous definition $1/s$ is the step of the discretization of the time variable, hence 
$(\eu^{s}(k),K^{s}(k))$ represents the approximate pair displacement-crack at the time $k/s$. 
We name any function $(\eu^{s}, K^{s}): N \rightarrow M$ an incremental solution if it 
verifies i) and ii) from the definition 4.1.. 

When $s_{i}$ converges to $\infty$ the time step goes to 0 and the incremental   
solution  $(\eu^{s_{i}},K^{s_{i}})([s_{i}t])$ converges to $(\eu,K)(t)$, for any $t > 0$.

A necessary condition for the existence of an energy minimizing movement introduced in 
definition 4.1. is that for any given $s$ the incremental solution 
$k \in N \  \mapsto \ (\eu^{s}, K^{s})(k)$ exists. The following theorem provides an answer to this 
existence query. In the general case $n=3$ this theorem is true, to our knowledge, 
only in a weak form, presented  in section 6. In the anti-plane case, however, due to  partial 
regularity results for the minimizers of Mumford-Shah functional from [DGCL], the theorem is
certainly true.   

\vspace{.5cm}

{\bf Theorem 4.1.} {\it Let $\Omega \subset\mathbb{R}^{n}$ be a bounded open set with piecewise smooth
boundary, let $(\eu_{0},K)$ be a given admissible pair displacement-crack in $\Omega$ and let 
$$\eu_{0}: N \ \rightarrow H^{\frac{1}{2}}(\paromega,R^{n}) \cap L^{\infty}(\paromega,R^{n})$$ 
be a given sequence of imposed
displacements such that $\eu_{0} \ = \ \eu(\eu_{0}(0),K)$ on $\paromega \setminus K$. 

Then there exists the sequence $(\eu,K): N \ \rightarrow \ M$ such that: 

\hspace{1.cm} i) $\eu(0) \ = \ \eu_{0}$ and $K(0) \ = \ K$; 

\hspace{1.cm} ii)  for any $k \in N$ there is a crack set $L(k+1)$ such that 
$(\eu(k+1),L(k+1)) \in M$, $\eu(k+1) \ = \ \eu_{0}(k+1)$ on $\paromega \setminus L(k+1)$ and 
$(\eu(k+1),L(k+1))$ is a minimizer of the functional 
$$(\ev, L) \in M \ , \ev \ = \ \eu_{0}(k+1) \mbox{ on } \paromega \setminus L \ \ \mapsto 
J\left( (\eu(k),K(k)), \ev,L) \right) \ \ .$$ 
The set $K(k+1)$ is given by the formula}  $$K(k+1) \ = \ K(k) \ \cup \ L(k+1) \ \ .$$ 

\vspace{.5cm}

\subsection{Features of the first model}
\indent

We shall investigate further the behaviour of the model proposed in definition 4.1. in the
particular case of anti-plane displacement. There are some obvious adjustments to be made. 
$\Omega$ is now a bounded domain in $R^{2}$ and the displacement is a scalar function $u$. 
The functional $J$ will take the form: 
$$J\left( (u,K), (v,L) \right) \ = \ \int_{\Omega} \mu \ \mid \nabla  v \mid^{2} \mbox{ d}x \ 
+ \ G \mathcal{ H}^{1}(L \setminus K) \ \ .$$ 

Let us consider a particular type of imposed displacement on $\paromega$. We split the boundary 
of the body in three parts: 
 $$\paromega = \overline{\gau^{1}} \cup \overline{\gau^{2}} \cup \overline{\gaf}$$
 $$\gau^{i} \cap \gaf = \emptyset \mbox{ , }  \overline{\gau^{1}} \cap \overline{\gau^{2}} = \emptyset 
\mbox{ , } \mathcal{ H}^{1}(\gau^{1}) \cdot \mathcal{ H}^{1}(\gau^{2}) \cdot \mathcal{ H}^{1}(\gaf) > 0\ \ .$$ 
At any moment $t \geq 0$, $\gaf$ is force free, i.e. the displacement is not prescribed on 
this part of the boundary. On $\gau^{1}$ and $\gau^{2}$ the imposed displacement is defined by 
the formula: 
$$u_{0}(t)(x) \ = \ \left\{ \begin{array}{ll}
0 & \mbox{ on } \gau^{1} \\
t \delta & \mbox{ on } \gau^{2} \ \ ,  
\end{array}
\right. $$
where $\delta$ is a positive constant with dimension of speed. This  displacement  is homogeneous 
with respect to the time variable: 
$$\forall t \ > \ 0 \ , \ u_{0}(t) \ = \ t u_{0}(1) \ \ .$$
We suppose further that at the moment $t=0$ there are no cracks in the body. This assumption 
takes the form $K \ = \ \emptyset$.   At $t=0$ we have $u_{0}(0) = 0$, hence the initial data are 
$(u_{0}=0, K=\emptyset)$. 

Let us consider a time discretization given by the parameter $1/s$ and the incremental solution 
$k \in N \  \mapsto \ (u^{s}, K^{s})(k)$ introduced in definition 4.1. for the initial data 
 and the imposed boundary described above. In order to shorten the notations we shall omit for 
the moment the superscript $s$.  

The incremental solution  $(u,K): N \rightarrow M$  is recursively  defined by the following two
rules: 

\hspace{1.cm} i) $u(0) \ = \ 0$ and $K(0) \ = \ \emptyset$; 

\hspace{1.cm} ii)  for any $k \in N$ we seek for the crack set $L(k+1)$ and for the displacement 
$u(k+1)$ such that 
$(u(k+1),L(k+1)) \in M$, $u(k+1) \ = \ (k+1)/s \ u_{0}(1)$ on 
$\left( \gau^{1} \cup \gau^{2} \right)  \setminus L(k+1)$ and 
$(u(k+1),L(k+1))$ is a minimizer of the functional 
$$(v, L) \in M \ , v \ = \ (k+1)/s \ u_{0}(1) \mbox{ on } 
\left( \gau^{1} \cup \gau^{2} \right)  \setminus L \ \ \mapsto 
J\left( (u(k),K(k)), (v,L) \right) \ \ .$$ 
The set $K(k+1)$ is given by the formula  $$K(k+1) \ = \ K(k) \ \cup \ L(k+1) \ \ .$$ 

Let us denote by $u_{\emptyset}$ the displacement of the body $\Omega$, without cracks,
 under the prescribed 
displacement on the boundary $u_{0}(1)$. With the use of a notation made before, 
$u_{\emptyset}$ is defined by $u_{\emptyset} \ = \ u(u_{0}(1),\emptyset)$. For 
any $k \in N$  we   have $(k/s \ u_{\emptyset} ,\emptyset) \in M$ and $k/s \ u_{\emptyset} \ = \ 
k/s \ u_{0}(1)$ on  $\gau^{1} \cup \gau^{2}$. Therefore, for any $k \in N$ 
we have $$J\left( (u(k),K(k)), (u(k+1),L(k+1)) \right) \ \leq \ 
J\left( (u(k),K(k)), ((k+1)/s \ u_{\emptyset},\emptyset) \right) \ \ .$$ The last inequality reads: 
\begin{equation}
\int_{\Omega} \mu \ \mid \nabla u(k+1) \mid^{2} \mbox{ d}x  \  + \ G \mathcal{ H}^{1} 
(L(k+1) \setminus K(k) ) \ \leq \ 
\left(\frac{k}{s}\right)^{2} \ \int_{\Omega} \mu \ \mid \nabla u_{\emptyset} \mid^{2} 
\mbox{ d} x  \ \ .
\label{ineq1}
\end{equation}

We can always find a curve in $\overline{\Omega}$ which separates $\gau^{1}$ from $\gau^{2}$. 
Moreover, we can find such a curve which is a length minimizer in the family of all curves 
 in $\overline{\Omega}$ separating $\gau^{1}$ from $\gau^{2}$. Let us denote this curve  by $S$
(which exists but it might not be unique). The domain $\overline{\Omega}$ has the following
decomposition with respect to $S$: 
$$\overline{\Omega} \ = \ \Omega^{1} \cup \Omega^{2} \ ,  \ \gau^{1} \subset \Omega^{1} \ , \ 
\gau^{2} \subset \Omega^{2} \ , \ \Omega^{1} \cap \Omega^{2} \ = \ \emptyset \ , $$
$$ \overline{\Omega^{1}} \cap \overline{\Omega^{2}} \ = \ S \ \ .$$
Let us define the following displacement: 
$$u_{S}(x) \ = \  \left\{ \begin{array}{ll}
0 & x \in \Omega^{1} \\
\delta & x \in \Omega^{2} \ \ \
\end{array}
\right. $$
It is easy to see that for any $k \in N$ we have $(k/s \ u_{S},S) \in M$ and 
$k/s \ u_{S} \ = \ k/s \ u_{0}(1)$ on $\left( \gau^{1} \cup \gau^{2} \right) \setminus S$. Therefore 
we obtain the following inequality: 
 \begin{equation}
\int_{\Omega} \mu \ \mid \nabla u(k+1) \mid^{2} \mbox{ d}x  \  + \ G \mathcal{ H}^{1} 
(L(k+1) \setminus K(k) ) \ \leq \ G \mathcal{ H}^{1}(S \setminus K(k)) \ \ .
\label{ineq2}
\end{equation}

From  (\ref{ineq2}) we derive the following conclusion: {\it for large time 
$k/s$ the crack set $K(k)$ is not void.} Indeed, suppose that the function 
$k \in N \ \mapsto \ (k/s \ u_{\emptyset}, \emptyset)$ is an incremental solution constructed 
by the rules i) and ii) above. Then for any $k \in N$ the inequality (\ref{ineq1}) becomes an 
equality and the inequality (\ref{ineq2}) takes the following form: 
\begin{equation}
\left(\frac{k}{s}\right)^{2} \ \int_{\Omega} \mu \ \mid \nabla u_{\emptyset} \mid^{2} 
\mbox{ d} x  \ \leq \ G \mathcal{ H}^{1}(S \setminus K(k)) \ \ ,
\label{ineq3}
\end{equation} 
which lead to contradiction. 
Therefore this model can predict crack appearance. The 
critical step  $k$, after which a crack appears in the body 
(as the incremental solution predict), is the greatest natural with the property (\ref{ineq3}).

The following theorem  contains stronger informations regarding the minimizers of the 
Mumford-Shah functional in our particular case.

\vspace{.5cm}

{\bf Theorem 4.2.} {\it Let $\Omega \subset\mathbb{R}^{2}$ be a bounded open set with piecewise smooth 
boundary $\paromega$ and let $\en$ be the field of outward normals over the boundary. 
Let us suppose that the boundary of $\Omega$ has the 
following decomposition: 
$$\paromega = \overline{\gau^{1}} \cup \overline{\gau^{2}} \cup \overline{\gaf}$$
 $$\gau^{i} \cap \gaf = \emptyset \mbox{ , }  \overline{\gau^{1}} \cap \overline{\gau^{2}} = \emptyset 
\mbox{ , } \mathcal{ H}^{1}(\gau^{1}) \cdot \mathcal{ H}^{1}(\gau^{2}) \cdot \mathcal{ H}^{1}\gaf) > 0\ \ .$$
Let us consider the functional 
$$I(v, K) \ = \ \frac{1}{2} \ \int_{\Omega} \mid \nabla v \mid \mbox{ d}x \ + \ G \mathcal{ H}^{1}(K 
) \ \ ,$$
defined over the set 
$$\left\{ (v,K) \mbox{ : } v \in C^{1}(\overline{\Omega} \setminus K,\mathbb{R}) \right\} \ \ .$$ 
Let, for any $D \in R$,  $u(D): \overline{\Omega} \rightarrow R$   
be  the solution of the problem: 
$$\left\{ \begin{array}{ll}
div \ \nabla v \ = \ 0 & \mbox{ in } \Omega \\
\nabla v  \en \ = \  0 & \mbox{ on } \gaf \\
v \ = \ 0 & \mbox{ on } \gau^{1} \\
v \ = \ D & \mbox{ on } \gau^{2} \ \ .  
\end{array} \right.$$
We suppose that exist strictly  positive numbers $c$ and $C$ such that for any $x \in \gau^{1} 
\cup \gau^{2}$ 
$$ C \geq \mid \nabla u(1) \en \mid (x) \geq c \ \ .$$

There exist then two numbers $m \leq M$, which depends only on $\Omega, \gau^{1},\gau^{2}$ and 
$\gaf$, such that: 

i) if $D^{2} < m$ then $(u(D),\emptyset)$ is the only minimizer of the functional $I$ over the set  
$$ M(D) \ = \ \left\{ (v,K) \mbox{ : } v \in C^{1}(\overline{\Omega} \setminus K,\mathbb{R}) \ , \ v = u(D) \mbox { on
} \gau^{1} \cup \gau^{2} \right\} \ \ , $$

ii) if $ D^{2} > M$ then any minimizer of the functional $I$ over the set $M(D)$ has the form 
$(u_{K}, K)$, with $\mid \nabla u_{K} \mid = 0$ almost everywhere in $\Omega$ and $K$ geodesic in $\Omega$ 
(i.e. length minimizer) separating $\gau^{1}$ from $\gau^{2}$. 

Moreover, if $c = C$ then $M=m$, hence if  $\nabla u(1) \en$ is piecewise
constant on $\gau^{1} \cup \gau^{2}$ then we have only two kinds of minimizing crack sets. }

\vspace{.5cm}

The theorem assures us that for small time $k/s$ the body remains sane and for large time $k/s$
a crack with a particular shape  appears in the body. Precisely, for small $k/s$ we have 
$(u(k),K(k)) \ = \ (k/s \ u_{\emptyset}, \emptyset)$ and for large $k/s$ we have 
$(u(k),K(k)) \ = \ (k/s \ u_{S},S)$. The theorem help us to find  particular cases 
when the passage from the first type of minimizer to the second one is brutal. Indeed, consider 
that $\Omega$ is a rectangle $(0,a) \times (0,L)$,  
$$\gau^{1} \ = \ (0,a) \times \left\{ 0 \right\} \ , \ 
\gau^{1} \ = \ (0,a) \times \left\{ L \right\} \ $$
and $\gaf$ is the remaining part of the boundary. Let us consider, for simplicity, that $\delta =
1$. With the notations from the theorem 4.2. we have $c\ = \ C$ therefore we have only two kinds 
of pairs displacement-crack  which compete. We use (\ref{ineq3}) in order to find the critical 
time $k/s$ when the incremental solution $(u,K):N \rightarrow M$ switches from  
$(k/s \ u_{\emptyset}, \emptyset)$ to $(k/s \ u_{S},S)$, where $S$ is, for example, 
$(0,a) \times \left\{ L/2 \right\}$. We find that 
 the critical $k/s$  is determined by the double inequality: 
$$\left( \frac{k}{s} \right)^{2} \ \leq \  \frac{G L}{\mu} \ \leq \ \left( \frac{k+1}{s}\right)^{2} 
\ \ .$$

We are lead to the definition of  the  critical moment $t_{c}$, given by the formula   
$$t_{c}^{2} \ = \  \frac{G L}{\mu} \ \ .$$  
$t_{c}$ is proportional with the square root of $L$. The anti-plane stress existing in the sane body 
at the moment $t_{c}$ has the following expression: 
$$\mu \nabla \left( t_{c} u_{\emptyset}  \right) \ = \ \left( 0, \frac{t_{c}}{L} \right) \ \ .$$
We can see that this stress depends on $L$, hence on the geometry of the body.  
Because $G$ is supposed to be a material constant we obtain  the following 
conclusion: {\it the model described above is not compatible with any model of crack appearance 
based on a critical stress as  material constant.}

\section{The improved model}
\indent

We have seen that the first model allows crack appearance but it is not compatible with 
any critical fracture stress based model. Our purpose is to improve the first model in 
order to allow the existence of a critical stress which damages a structure. We shall find 
a way to make this improvement by studying first how smooth brittle propagation of 
cracks can be described with the Mumford-Shah energetic functional.

\subsection{Smooth brittle crack propagation}

\indent

There are two steps in order to define the  notion of smooth brittle crack propagation. 
The first step consists in smoothness demands on the initial crack set $K$. We shall 
suppose that $K$ is endowed with the structure of  manifold with boundary. The boundary of 
$K$, denoted by $\partial K$, represents the edge of the crack. The second step consists in 
smoothness demands on the evolution $t \mapsto K_{t}$ of the crack. We shall restrict our attention 
only to evolutions of the initial crack $K$ obtained by smooth deformations of $K$. The initial 
crack may be as complex as we wish, because the structure of  manifold with boundary allows that, 
but this complexity remains the same during the propagation of the crack.

We shall work with  deformations of the initial crack set $K$ by endomorphisms of $\Omega$. Let 
us consider the following set of diffeomorphisms: 
\begin{equation}
\mathcal{ D}^{s} \ = \ \left\{ \phi \in C^{\infty}(\Omega,\Omega) \cap W^{s,2}(\Omega,R^{n}) 
\mbox{ : } \phi^{-1} \in C^{\infty}(\Omega,\Omega) \mbox{ and } supp \ (\phi - 1_{\Omega}) \subset 
\Omega \right\} \ \ .
\end{equation}  
We have denoted by $n$ the dimension of the space where $\Omega$ lies. The introduction of 
the Sobolev space $W^{s,2}(\Omega,R^{n})$ has been made for mathematical reasons ( to be found 
for example in Ebin \& Marsden [EbM]) and the number $s$ 
is chosen to be greater than $\frac{n}{2}+2$. In the paper [Bu3] a rigorous 
mathematical description of smooth brittle crack propagation can be found. We mention 
that the number $s$ controls the variation of the smoothness of the deformed crack set $\phi(k)$ 
from the smoothness of the initial crack set $K$. 

The condition $ supp \ (\phi - 1_{\Omega}) \subset \Omega$ means that near the boundary of 
$\Omega$ $\phi$ equals the identity map. 

\vspace{.5cm}

{\bf Definition 5.1.} {\it A smooth fracture curve is a function 
$$t \in [0,T] \  \mapsto \ \phi_{t} \in  \mathcal{ D}^{s} \ \ ,$$
which has the following properties: 

\hspace{2.cm} i) $\phi_{0} \ = \ 1_{\Omega}$ ,

\hspace{2.cm} ii) the map $t \in [0,T] \  \mapsto \ \phi_{t} \in  \mathcal{ D}^{s}$ is continuous 
with respect to the topology induced by the norm $max \ (\| \cdot \|_{L^{\infty}}, \| \cdot 
\|_{W^{s,2}} )$; for every $t \in [0,T]$ $\dot{\phi_{t}}$ exists and 
$$\eta_{t} \ = \ \dot{\phi_{t}}. \phi_{t}^{-1} \ \in \ W^{s,2}(\Omega,R^{n}) \cap
L^{\infty}(\Omega,R^{n}) \ \ ,$$

\hspace{2.cm} iii) for any $t < t'$ we have $\phi_{t}(K) \subset \phi_{t'}(K)$ .

We have used the notation $f.g$ for the composition of the function $f$ with $g$. 

A crack evolution curve is associated to the smooth fracture curve $t \mapsto \phi_{t}$ and 
initial crack $K$ by the formula: }
$$K_{t} \ = \ \phi_{t}(K) \ \ .$$

\vspace{.5cm}

There are infinitely many smooth crack propagation curves $t \mapsto \phi_{t}$ 
with the same associated crack evolution curve $t \mapsto K_{t}$. 

Under smoothness assumptions on the initial crack set $K$, 
for a smooth crack propagation curve $t \mapsto \phi_{t}$ the condition that the crack 
grows implies that for any $t \geq 0$ we have: 
\begin{equation}
\left\{ \begin{array}{ll}
\eta_{t} \cdot \en \ = 0 & \mbox{ on } \phi_{t}(K) \\
\int_{\phi_{t}(K)} div_{s} \eta_{t} \mbox{ d}\hen \ \geq 0 \ \ .
\end{array} \right.
\label{difcond}
\end{equation}
The integral from the last inequality equals the variation of  the  area of the crack set (see 
Allard [All]). 
The operator 
$div_{s}$ is the tangential derivative with respect to the surface $\phi_{t}(K)$ and it has the 
following form: 
$$div_{s} \eta \ = \ div \ \eta \ - \en \cdot \left( \nabla \eta \right) \en \ \ ,$$ 
where $\en$ is the normal to the surface $\phi_{t}(K)$.

\subsection{K2 functional and the Griffith criterion} 

\indent

We want now to reformulate the Griffith criterion of brittle crack propagation 
(\ref{grif}) in terms of smooth crack propagation curves. Our assumptions on the 
evolution of  the body are the following: 

\hspace{2.cm} A1) the evolution of the linear elastic body $\Omega$ is quasi-static, 

{\hspace{2.cm} A2) a smooth curve of imposed displacements $t \mapsto \eu_{0}(t)$ is given 
on the boundary $\paromega$. 

The assumption A2) can be modified by the replacement of $\paromega$ with a fixed part of 
the boundary $\gau$; on the remaining part $\gaf$ we suppose  that the body is force free. 

At any moment $t$ the state of the body is described by the admissible  pair displacement-crack 
$(\eu(t),K_{t})$. The assumption A1) implies that the displacement $\eu(t)$ is determined by 
the knowledge of $K_{t}$ and boundary condition $\eu_{0}(t)$.  With a notation used several times 
before, we have $\eu(t) \ = \ \eu ( \eu_{0}(t), \phi_{t}(K))$. We shall change for our purposes 
this notation by writing: 
$$\eu(t) \ = \ \eu ( \eu_{0}(t), \phi_{t}) \ \ .$$ 

\vspace{.5cm}

{\bf Definition 5.2.} {\it Let $t \in [0,T] \mapsto \eu^{0}(t) \in C(\paromega,R^{n})$ be a 
$C^{1}$ curve of imposed
displacements on the exterior boundary of the body. A balanced fracture curve is any 
$C^{1}$ function 
$$t \in [0,T] \mapsto (\eu^{*}_{t},\phi_{t}) \in W^{1,2}(\Omega \setminus K,R^{n}) 
 \times \mathcal{ D}^{s} \
\ ,$$ 
satisfying the following items:

\hspace{2.cm} i) for any $t \in [0,T]$  we have 
$$\eu^{*}_{t} . \phi_{t}^{-1} \ = \ \eu(t) = \eu ( \eu_{0}(t), \phi_{t}) \ \ ,$$

\hspace{2.cm} ii) $t \mapsto \phi_{t}$ is a smooth crack propagation curve.}

\vspace{.5cm}

For given curve of boundary displacement and initial crack set $K$, for any smooth crack 
propagation curve there is only one associated balanced fracture curve. The Griffith criterion 
of brittle crack propagation will act as a selection criterion amongst all smooth crack propagation
curves.  

We have seen that an equivalent form of the Griffith criterion  is 
(\ref{grifnew}). With the change of notation $\eT(\phi_{t}) \ = \ \eT(\phi_{t}(K))$, we say that 
a smooth crack propagation curve is compatible with the Griffith criterion if for any  $t$ we have 
\begin{equation}
 \frac{1}{2}  \langle \frac{d}{dt} \left[ \eT(\phi_{t}) \right] \eu_{0}(t),\eu_{0}(t) \rangle \ + \
G \frac{d}{dt}\left\{
\hen(\phi_{t}(K)) \right\} \ \leq \ 0 \ \ .
\label{grifphi}
\end{equation}

The first term from (\ref{grifphi}) represents the variation of the elastic energy of the body 
calculated for the following variation of the displacement: 
$$\tau \ \mapsto \ \eu(t+\tau) \ = \ \eu(\eu_{0}(t),\phi_{t+\tau}) \ \ .$$ 
The dependence of $\eu(t+\tau)$ with respect to $\phi(t+\tau)$ is implicit. An explicit variation 
of the displacement  
would be preferable, like this one: 
$$\ev(\tau) \ = \ \eu(t).\phi_{t} . \phi_{t+\tau}^{-1} \ \ .$$
We are lead, by a change of variables, to the following equality: 
\begin{equation}
2 \ \frac{d}{d\tau} \ \left( \int_{\Omega} w(\nabla \ev(\tau)) \mbox{ d}x \right)_{|_{\tau=0}} 
\ = \  
\int_{\Omega} \left\{ \left[ \eC \nabla \eu(\eu_{0}(t),\phi_{t}) : \nabla \eu(\eu_{0}(t),
\phi_{t}) \right] div \ \eta_{t} \ - \right. 
\label{apk2}
\end{equation}
$$ \left. - 2 \left[ \eC \nabla \eu(\eu_{0}(t),\phi_{t}) \right]_{ij} 
\left[ \nabla \eu(\eu_{0}(t),\phi_{t}) \right]_{ik} \left[ \nabla \eta_{t} \right]_{kj} 
\right\} \mbox{ d}x \ \ .$$

\vspace{.5cm}

{\bf Proposition 5.1.} {\it Let $t \mapsto \phi_{t}$ be a smooth crack propagation curve, 
$\eta_{t} \ = \ \dot{ \phi_{t}} . \phi_{t}^{-1}$ and $\eu_{0} \in H^{\frac{1}{2}}(\paromega) 
\cap L^{\infty}(\paromega)$. Let us define, for fixed $t$,   
$\ev(\tau) \ = \ \eu(t).\phi_{t} . \phi_{t+\tau}^{-1}$ . Then  the following inequality is true : }
\begin{equation}
 \langle \frac{d}{dt} \left[  \eT(\phi_{t}) \right] \eu_{0} , \eu_{0}  \rangle \ \leq  \ 2 \ 
\frac{d}{d\tau} \ \int_{\Omega} w(\nabla \ev(\tau)) \mbox{ d}x_{|_{\tau=0}}  \ \ .
\label{p1}
\end{equation}

\vspace{.5cm}

This proposition, together with the equality (\ref{apk2}),
 allows us to introduce a generalization of the J integral. We use  the notation 
$\eta \in W^{s,2}_{0}(\Omega,R^{n})$ for $\eta \in  W^{s,2}_{0}(\Omega,R^{n})$ with  
null trace on  $\paromega$. 

\vspace{.5cm}

{\bf Definition 5.3.}{ \it The generalized  J integral is the  following functional }
$$K2 \ : \ \mathcal{ D}^{s} \times \left\{ W^{s,2}_{0}(\Omega,R^{n}) \cap L^{\infty}(\Omega,R^{n}) 
\right\}  
\times \left\{  H^{\frac{1}{2}}(\paromega) \ \cap L^{\infty}(\paromega,R^{n}) 
\right\} \ \rightarrow \\mathbb{R}
\ \ ,$$
\begin{equation}
K2(\phi,\eta,\eu_{0}) \ = \ 
\int_{\Omega} \left\{ - \ \frac{1}{2} \left[ \eC \nabla \eu(\eu_{0},\phi) : \nabla \eu(\eu_{0}
,\phi) \right] div \ \eta \ + \right.
\label{k2}
\end{equation}
$$\left.  +  \left[ \eC \nabla \eu(\eu_{0},\phi) \right]_{ij} 
\left[ \nabla \eu(\eu_{0},\phi) \right]_{ik} \left[ \nabla \eta \right]_{kj} 
\right\} \mbox{ d}x \ \ .$$  

\vspace{.5cm}

In order to explain why the functional $K2$ is the generalization J integral, we begin by 
a temporary introduction of easier notations: 
$$\sig = \eC \nabla \eu(\eu_{0},\phi) \ \ , \  \eu = \eu(\eu_{0},\phi) \ \ , \ w = \frac{1}{2} 
\eC \nabla \eu(\eu_{0},\phi) : \nabla \eu(\eu_{0},\phi) \ \ .$$ 

We define the tubular neighbourhood, of radius $r$, 
 of the edge $\partial \phi(K)$ of the crack set $K$: 
$$B_{r}= B_{r}(\partial \phi(K)) = \cup_{x \in \partial \phi(K)} B(x,r) \ \ .$$
The field of normals over $\partial \ B_{r}(\partial \phi(K))$ will be denoted by $\nu$, without 
specifying the parameter $r$. 

If $\eu$ belongs to  $C^{2}$ then  we have
\begin{equation}
w \ \eta_{i,i} - \sigma_{ij} \eu_{i,k} \ \eta_{k,j} = 
\left[ w \ \eta_{i} \ - \ \sigma_{lj} \eu_{l,k} \
\eta_{k} \right]_{,i} - \sigma_{km} \eu_{k,mi} \ \eta_{i} + \sigma_{li} \eu_{l,ki}  \ \eta_{k} + 
\sigma_{li,i} \eu_{l,k}  \ \eta_{k} \ \ .
\label{nece}
\end{equation}

According to the assumption A1), the divergence of the stress field $\sig$ equals $0$. We
integrate the  equality (\ref{nece}) over $\Omega \setminus B_{r}$ and we obtain: 
$$K2(\phi,\eta,\eu_{0}) \ = \ - \ \lim_{r \rightarrow 0} \ 
\int_{\Omega \setminus B_{r}} \left[w \ \eta_{i} \ - \ \sigma_{li} \eu_{l,k} \
\eta_{k} \right]_{,i} \mbox{ d}x \ \ .$$ 
By a flux-divergence formula, we are lead to the following expression of $K2$: 
$$K2(\phi, \eta,\eu_{0}) = \lim_{r \rightarrow 0} \left\{ \int_{\partial B_{r}(\partial \phi(K))} 
\left\{ - w \eta \cdot \nu \ + \sigma_{li} \eu_{l,k} \ \eta_{k} \ \nu_{i} \right\} \right\} + \int_{\phi(K)} 
[w] \eta \cdot \en \mbox{ d} \hen  \ \ .$$
The functional $K2$ is interesting in the case when $\eta \cdot \en = 0$ on $\phi(K)$, as 
(\ref{difcond}) suggests. In this case we have: 
\begin{equation}
 K2(\phi, \eta,\eu_{0}) = \lim_{r \rightarrow 0} \left\{ \int_{\partial B_{r}(\partial \phi(K))} 
 - \ \left\{ w \eta \cdot \nu \ - \sigma_{li} \eu_{l,k} \ \eta_{k} \ \nu_{i} \right\} \mbox{ d} \hen  
\right\} \ \ . 
\label{k2asrice}
\end{equation}
Let us consider that we are in the case of plane displacements (hence $n=2$) and that the crack 
set $\phi(K)$ lies on the $Ox_{1}$ axis. If we take $\eta$  equal to $(1,0)$ in a
neighbourhood of the edge of the crack then we have:
$$ K2(\phi, \eta,\eu_{0}) = \lim_{r \rightarrow 0} \left\{ \int_{\partial B_{r}(\partial \phi(K))} 
 - \ \left\{ w \nu_{1} \ - \sigma_{ki} \eu_{k,1} \  \ \nu_{i} \right\} \mbox{ d} \hen  
\right\} \ \ . $$
We recognize in the right term of the equality above the expression of the classical J integral.

We propose the following selection criterion for smooth crack propagation curves: 

\vspace{.5cm}

\hspace{2.cm} {\it A smooth crack propagation curve $t \mapsto \phi_{t}$ satisfies the generalized 
Griffith criterion if at any moment $t \geq 0$} we have $\eta_{t} \cdot \en \ = \ 0$ on 
$\phi_{t}(K)$ and 
\begin{equation}
K2(\phi_{t},\eta_{t},\eu_{0}(t)) \ \geq \ G \ \frac{d}{dt} \hen(\phi_{t}(K)) \ \ .
\label{g2}
\end{equation}

\vspace{.5cm}

Ohtsuka [Oht1---4] proves that under stronger smoothness assumptions on $K$ and on the curve 
$t \mapsto K_{t}$ always exists a smooth crack propagation curve $t \mapsto \phi_{t}$ such 
that for any $t$ and with our notations  we have:  
$$\langle \frac{d}{dt} \left[  \eT(\phi_{t}) \right] \eu_{0}(t) , \eu_{0}(t) \rangle \ + \ 
K2(\phi_{t},\eta_{t},\eu_{0}(t)) \ = \ 0 \ \ .$$
For this reason we consider that (\ref{g2}) is not too strong with respect to the classical Griffith 
criterion.

\subsection{Extension of K2 and admissible cracks}

\indent

We want to extend the Griffith criterion of brittle fracture propagation (\ref{g2}) in order to 
allow crack appearance. The leading idea is to consider crack evolution curves $t \mapsto K_{t}$ 
which are limits of crack evolution curves of the form $t \mapsto \phi_{t}(K)$. 

Let $t \mapsto \phi_{t}$ be a smooth crack propagation curve. At any moment $t$ the vector field 
$\eta_{t}=\dot{\phi_{t}}.\phi_{t}^{-1}$ represents the propagation speed of the edge $\partial 
\phi_{t}(K)$ of the crack $\phi_{t}(K)$. Precisely the restriction  of $\eta_{t}$ to 
$\partial 
\phi_{t}(K)$ represents the distribution of  speed of propagation of the points 
belonging to this $n-2$ surface. 
The appearance of a new crack at the moment $t$ 
is seen as a limit of  processes of smooth crack  propagation, when the distribution of speed 
$\eta_{t}$ develops jumps. 

We shall consider therefore a sequence $(\eta_{h})_{h}$ in the space $W^{s,2}(\Omega,R^{n}) 
\cap L^{\infty}(\Omega,R^{n})$. For each $h$ we define the following flow 
$\tau \mapsto \phi^{h}_{\tau} $:
$$\phi^{h}_{0} \ = \ 1_{\Omega} \ , \ \ \phi^{h}_{\tau} = \ 1_{\Omega} \ + \ \tau \eta_{h} \ \ .$$
For small times $\tau$ we have $\phi^{h}_{\tau} \in \mathcal{ D}^{s}$ and for any $h$ we see that 
$$\dot{\phi^{h}_{0}}.(\phi^{h}_{0})^{-1} \ = \ \eta_{h} \ \ .$$ 
Let us suppose that $\eta_{h}$ converges almost everywhere to $\eta$. Then for any $\tau$ 
$\phi^{h}_{\tau}$ converges almost everywhere to $\phi_{\tau} \ = \  1_{\Omega} \ + \ \tau \eta$. 
We make the following supplementary assumptions:

\hspace{2.cm} S1) for small $\tau$ $\phi_{\tau}$ is almost everywhere injective,

\hspace{2.cm} S2)  $\Omega \setminus \phi_{\tau}(\Omega)$ has finite Hausdorff $n-1$ measure,

\hspace{2.cm} S2) let us denote by $S_{\eta}$ the set where $\eta$ has no approximate limit 
(the complementary of the Lebesgue set of $\eta$, see for this  
the section of proofs) and by $\mid Df \mid $ the total variation measure associated to the 
distributional derivative  of the 
function $f$; then we have
$$\mid D^{s}\eta \mid (\Omega \setminus   S_{\eta}) \ = \ 0 \ \ .$$
The assumption S1) assures us that $\phi_{\tau}$ almost everywhere maps different  points from 
$\Omega$ 
in different places in $\phi_{\tau}(\Omega)$; 
S2) prevents the case where in the limit appear holes in the configuration $\phi_{\tau}(\Omega)$
  and 
S3) is a more sophisticated condition which says that no strange Cantor sets appear in  
$\phi_{\tau}(\Omega)$. 

Under these assumptions it is easy to prove that on $S_{\eta}$ the jump of $\eta$ satisfies the 
relation: 
$$[\eta]\cdot \en \ = \ 0 \ \ .$$
Indeed, suppose that in a neighbourhood of $x \in S_{\eta}$ we have $[\eta]\cdot \en \ < \ 0$. Then 
the assumption S2) is contradicted because a solid neighbourhood of $x$ is transformed by 
$\phi_{\tau}$ in a neighbourhood with a hole, when $\tau > 0$; if 
we have  $[\eta]\cdot \en \ > \ 0$ then S1) is contradicted because even if local injectivity is
respected, the global injectivity in the form S1) is not.   

We shall consider therefore pairs $(\eta,N)$ where $N$ is a topologically closed countably 
rectifiable set,  $\eta$ has tangential jumps on the surface $N$ and  
satisfies the smoothness 
assumption $\eta \in W^{s,2}_{0}(\Omega \setminus N,R^{n}) \cap L^{\infty}(\Omega,R^{n})$. 
The subscript $0$ in the notation $W^{s,2}_{0}(\Omega \setminus N,R^{n})$ means that $\eta = 0$ 
on $\paromega$ in the sense of traces. 

Let us suppose that  there is no initial crack in the body: $K \ = \ \emptyset$. We perform the
same calculation for $K2(1_{\Omega},\eta,\eu_{0})$ as we did after definition 5.3. and we obtain 
from (\ref{nece}) the expression: 
$$K2(1_{\Omega},\eta,\eu_{0}) \ = \ - \ 
\int_{\Omega} \left[ w \ \eta_{i} \ - \ \sigma_{lj} \eu_{l,k} \
\eta_{k} \right]_{,i} \ \ .$$
Let us suppose that $N$ is a surface with boundary and  let $B_{r}$ be a tubular neighbourhood  
of $\partial N$, of radius $r$. Because of the assumption  $[\eta]\cdot \en \ = \ 0$ on $N$, we
obtain: 
$$ K2(1_{\Omega},\eta,\eu_{0}) \ = \ \int_{N} \sigma_{li} \en_{i} \eu_{l,k} [\eta_{k}] \mbox{ 
d}\hen \ \ .$$
It is natural to try to modify the Griffith criterion (\ref{g2}) in order to have a control on the 
integral from above. We propose the following differential criterion of crack appearance (DA), which 
make a selection amongst all crack sets which can appear in the body. The constant $\Sigma$, 
with the dimension of a stress, which appears in this criterion, is postulated to be a 
constant of material.  

\vspace{.5cm}

{\bf (DA).} {\it Let us consider the elastic body $\Omega$ and the imposed boundary displacement 
$\eu_{0}$. A crack set $N$ can appear in the body if there exists a vector field 
$\eta \in W^{s,2}_{0}(\Omega \setminus N,R^{n}) \cap L^{\infty}(\Omega,R^{n})$, such that 
$[\eta]\cdot \en \ = \ 0$ on $N$ and} 
\begin{equation}
\int_{N} \sigma_{li}(\eu_{0},\emptyset) \en_{i} \eu_{l,k}(\eu_{0},\emptyset) [\eta_{k}] \mbox{ 
d}\hen  \ \geq \ \| \eta \|_{L^{\infty}} \ \Sigma \ \hen(N) \ \ .
\label{da}
\end{equation}

\vspace{.5cm}

We give also a criterion of local crack appearance (LA). This criterion tells us if in the point 
$x \in \Omega$ a crack with normal $\en$ can appear. 

\vspace{.5cm}

{\bf (LA).} {\it Let us consider the elastic body $\Omega$ and the imposed boundary displacement 
$\eu_{0}$. In the point $x \in \Omega$ a crack with normal $\en$ can appear if }
\begin{equation}
sup \ \left\{ \sigma_{li}(\eu_{0},\emptyset) \en_{i} \eu_{l,k}(\eu_{0},\emptyset) \nu_{k} 
\mbox{ : }  \nu \in\mathbb{R}^{n} \ , \ \mid \nu \mid = 1 \ , \ \en \cdot \nu = 0  \right\} \ \geq \ \Sigma \ \ .
\label{la}
\end{equation}

\vspace{.5cm}

From (\ref{da}) and (\ref{la}) we see that if $N$ is smooth enough and if 
 for any $x \in N$, the criterion (LA) is satisfied 
for the pair $(x,\en(x))$, where $\en(x)$ is the normal to $N$ at $x$, then $N$ satisfies the 
global criterion (DA). 

Let us suppose that the body $\Omega$ is  a cylinder $\omega \times [0,L]$ and $\eu_{0}$ imposed 
on the top and bottom of this cylinder such that 
$$\eu(\eu_{0},\emptyset) (x_{1},x_{2},x_{3}) \ = \ (0,0, a x_{3}) \ \ , a > 0 \ \ .$$
The stress $\sig(\eu_{0},\emptyset)$ has the form: 
$$\sig(\eu_{0},\emptyset) \ = \ \left( \begin{array}{ccc}
0 & 0 & 0 \\
0 & 0 & 0 \\
0 & 0 & \sigma 
\end{array} \right) \ \ ,$$
where $\sigma \ = \ Ea$. 
If we denote by $\alpha$ the angle between $\en$ and the $Ox_{3}$ axis, we have:  
$$ sup \ \left\{ \sigma_{li}(\eu_{0},\emptyset) \en_{i} \eu_{l,k}(\eu_{0},\emptyset) \nu_{i} 
\mbox{ : }  \nu \in\mathbb{R}^{n} \ , \ \mid \nu \mid = 1 \right\} \ = \   \frac{1}{2E}\sigma^{2} \ sin(2\alpha) 
\ \ .$$ 
The maximum value of this expression is attained for $\alpha=\pi / 4$. Therefore in the experience 
of uniaxial traction the (LA) criterion affirms that a crack can appear if
\begin{equation} 
\frac{1}{2E}\sigma^{2} \ \geq  \ \Sigma \ \ ,
\label{sic}
\end{equation}
and if we have equality in the relation above then the normal of the crack predicted by (LA) makes 
the angle $\pi / 4$ with the axis of the cylinder. The relation 
(\ref{sic}) gives us the value of the critical stress for uniaxial traction 
(which is a constant of material this time): 
\begin{equation}
\sigma_{cr} \ = \ \sqrt{2E\Sigma} \ \ .
\label{cris}
\end{equation}

\subsection{The improved model}
\indent 

In this section we propose an improved energy minimizing movement formulation to the problem 
of brittle fracture evolution (\ref{fracpro}). The model is based on two constants of material 
connected to fracture, namely $G$ and $\Sigma$ previously introduced. In this formulation the 
critical stress which lead to fracture in a traction experiment, defined by  (\ref{cris}),  
is a constant of material. 

We denote by $S_{n}$ the set of all $\nu \in\mathbb{R}^{n}$ with $\mid \nu \mid = 1$. 
Let us define the following function: 
$$f_{\infty} \ : \\mathbb{R}^{n \times n} \times\mathbb{R}^{n \times n} \times S_{n} \rightarrow\mathbb{R}\cup \left\{ + \infty 
\right\} \ \ ,$$
\begin{equation}
f_{\infty} (\sig,\eF, \en) \ = \ \left\{ \begin{array}{ll}
G & \mbox{ if } sup \ \left\{ \sigma_{li} \en_{i} \eF_{l,k} \nu_{k} 
\mbox{ : }  \nu \in\mathbb{R}^{n} \ , \ \mid \nu \mid = 1 \ , \ \en \cdot \nu = 0  \right\} \ \geq \  
\Sigma \\
+ \infty & \mbox{ otherwise } \ \ . 
\end{array}
\right. 
\label{finfty}
\end{equation}
The physical dimension of $f_{\infty}$ is the same as the one of $G$. 

With the use of the function $f_{\infty}$ the 
criterion of brittle crack appearance (LA), takes the form: {\it given the imposed boundary 
displacement $\eu_{0}$, $x \in \Omega$ and $\en \in S_{n}$, a crack of normal $\en$ can pass 
by $x$ if: }
$$f_{\infty}(\sig(\eu_{0},\emptyset),\nabla \eu (\eu^{0},\emptyset),\en) \ < \ + \infty \ \ .$$
The idea of the improved model is to consider only pairs displacement-crack $(\ev,L) \in M$
admissible with respect to (LA). 

\vspace{.5cm}

{\bf Definition 5.4.} {\it Let us consider the space M of admissible pairs displacement-crack 
endowed with the topology given by the convergence: 
$$(\eu_{h},K_{h}) \ \rightarrow \ (\eu,K) \mbox{ if } 
\left\{ \begin{array}{l} 
\eu_{h} \ L^{2} \rightarrow \ \eu \\
\hen(K_{h} \Delta K) \ \rightarrow \ 0 \ \ . 
\end{array} \right.$$
We define the functions 
$$J_{\infty}: M \times M \rightarrow\mathbb{R}\ ,$$
$$  J_{\infty}\left( (\eu,K),(\ev,L)\right) \ = \ 
\int_{\Omega} w(\nabla \ev ) \mbox{ d}x \ + \ \int_{L \setminus K}f_{\infty}(\sig(\eu), 
\nabla \eu, \en) \mbox{ d}\hen(x) \ \ ,$$
$$\Psi: [0,\infty) \times M \rightarrow \left\{0,+\infty \right\} \ , \ \ \Psi(\lambda, (\eu,K)) 
\ = \ \left\{ \begin{array}{ll} 
0 & \mbox{ if } \ev = \eu_{0}(\lambda) \mbox{ on } \paromega \setminus K \\
+\infty & \mbox{ otherwise } \ \ . 
\end{array} \right.$$
We consider the initial data $(\eu_{0},K) \in M$ such that $\eu_{0} \ = \ \eu(\eu_{0}(0),K)$. 

For any $s \geq 1$ we recursively define $(\eu^{s}, K^{s}): N \rightarrow M$ like this: 

\hspace{1.cm} i) $(\eu^{s}, K^{s})(0) \ = \ (\eu_{0},K)$ ;

\hspace{1.cm} ii) for any $k \in N$  $(\eu^{s}, L^{s})(k+1) \in M$  minimizes the functional 
$$ (\ev,L) \in M \ \mapsto \ J_{\infty}\left( ( (\eu^{s},K^{s})(k), (\ev,L) \right) \ + \ \Psi((k+1)/s,
(\ev,L))$$ 
over $M$.  $K^{s}(k+1)$ is defined by the formula: 
$$K^{s}(k+1) \ = \ K^{s}(k) \cup L^{s}(k+1) \ \ .$$

An energy minimizing movement associated to $J_{\infty}$ with the constraints (\ref{cevol}), $\Psi$ and 
initial data $(\eu_{0},K)$  is 
any $(\eu,K): [0,+\infty) \rightarrow M$ having the property:  
there is a diverging sequence $(s_{i})$ such that for any $t > 0$ 
$$(\eu^{s_{i}},K^{s_{i}})([s_{i}t]) \ \rightarrow (\eu,K)(t) \  \mbox{ as } i  \rightarrow \infty \
\ .$$}

\vspace{.5cm}

We are interested if for fixed $s$ an incremental solution exists. There is no result to our
knowledge that assures the existence of a minimizer of the functional 
$$(\ev,L) \in M \ \mapsto \ \int_{\Omega} w(\nabla \ev) \mbox{ d}x \ + \ 
\int_{L \setminus K}f_{\infty}(\sig(\eu), 
\nabla \eu, \en) \mbox{ d}\hen(x) \ \ $$ 
in our case. That is why we prefer to modify the function $f_{\infty}$. This function imposes 
a cost equal to $+\infty$ to the pairs displacement-crack which are not compatible with 
the (LA) criterion. We shall demand a finite but great cost instead, hoping that non admissible 
pairs will not enter in competition with admissible ones. 

Let us consider a number $C > G$ (with the same physical dimension as $G$) and a function 
$$f_{C} \ : \\mathbb{R}^{n \times n} \times\mathbb{R}^{n \times n} \times\mathbb{R}^{n} \rightarrow\mathbb{R} \ \ ,$$
with the following properties: 

\hspace{2.cm} i) $f_{C}$ is positively 1-homogeneous with respect to the third variable, 

\hspace{2.cm} ii) for any $\en \in\mathbb{R}^{n}$ such that $\mid \en \mid = 1$, if 
$$ sup \ \left\{ \sigma_{li} \en_{i} \eF_{l,k} \nu_{k} 
\mbox{ : }  \nu \in\mathbb{R}^{n} \ , \ \mid \nu \mid = 1 \ , \ \en \cdot \nu = 0  \right\} \ \geq \  
\Sigma $$
then $f_{C}(\sig,\eF,\en) = G$, 

\hspace{2.cm} iii) for any $\en \in\mathbb{R}^{n}$ such that $\mid \en \mid = 1$, if 
$$ sup \ \left\{ \sigma_{li} \en_{i} \eF_{l,k} \nu_{k} 
\mbox{ : }  \nu \in\mathbb{R}^{n} \ , \ \mid \nu \mid = 1 \ , \ \en \cdot \nu = 0  \right\} \ < \  
\Sigma $$
then $f_{C}(\sig,\eF,\en) > G$, 
 
\hspace{2.cm} iv) for any $(\sig,\eF,\en) \in \mathbb{R}^{n \times n} \times\mathbb{R}^{n \times n} \times S_{n}$ 
we have $f_{C}(\sig,\eF,\en) \leq C$ . 

As for theorem 4.1., the following statement has been proven to be true only in a weak sense,
described in the next section. 

\vspace{.5cm} 

{\bf Theorem 5.2.} {\it  Let $\Omega \subset\mathbb{R}^{n}$ be a bounded open set with piecewise smooth
boundary, let $(\eu_{0},K)$ be a given admissible pair displacement-crack in $\Omega$ and let 
$$\eu_{0}: N \ \rightarrow H^{\frac{1}{2}}(\paromega,R^{n}) \cap L^{\infty}(\paromega,R^{n})$$ 
be a given sequence of imposed
displacements such that $\eu_{0} \ = \ \eu(\eu_{0}(0),K)$ on $\paromega \setminus K$. 
 Let us consider  the functional $J_{C}$ 
$$J_{C}: M \times M \rightarrow\mathbb{R}\ ,$$
$$  J_{C}\left( (\eu,K),(\ev,L)\right) \ = \ 
\int_{\Omega} w(\nabla \ev ) \mbox{ d}x \ + \ \int_{L \setminus K}f_{C}(\sig(\eu), 
\nabla \eu, \en) \mbox{ d}\hen(x) \ \ ,$$
 where $f_{C}$ is chosen to satisfy the assumptions i)---iv) from above. 
Then there exists the sequence $(\eu,K): N \ \rightarrow \ M$ such that: 

\hspace{1.cm} i) $\eu(0) \ = \ \eu_{0}$ and $K(0) \ = \ K$; 

\hspace{1.cm} ii)  for any $k \in N$ there is a crack set $L(k+1)$ such that 
$(\eu(k+1),L(k+1)) \in M$, $\eu(k+1) \ = \ \eu_{0}(k+1)$ on $\paromega \setminus L(k+1)$ and 
$(\eu(k+1),L(k+1))$ is a minimizer of the functional 
$$(\ev, L) \in M \ , \ev \ = \ \eu_{0}(k+1) \mbox{ on } \paromega \setminus L \ \ \mapsto 
J_{C} \left( (\eu(k),K(k)), \ev,L) \right) \ \ .$$ 
The set $K(k+1)$ is given by the formula}  $$K(k+1) \ = \ K(k) \ \cup \ L(k+1) \ \ .$$ 

\vspace{.5cm}

\section{Proofs}
\indent

\subsection{Weak versions of theorems 4.1. and 5.2.}
\indent

This section is dedicated to a brief voyage trough the spaces $\eSBV$ and $\eSBD$. The 
weak forms of theorems 4.1. and 4.2. are direct applications of results listed below. 

The space $\eSBV(\Omega,R^{n})$ of special functions with bounded variation was introduced by 
De Giorgi and Ambrosio in the study of a class of free discontinuity problems ([DGA], [A1], [A2]). 
For any function $\eu \in L^{1}(\Omega,R^{n})$ let us denote by $D\eu$ the distributional 
derivative of $\eu$ seen as a vector measure. The variation of $D\eu$ is a scalar measure 
defined like this:  
for any Borel measurable subset $B$ of $\Omega$ the variation of $D\eu$ over $B$ is 
\begin{displaymath}
\mid D\eu \mid  (B) \ = \  sup \ \left\{ \sum^{\infty}_{i=1} \mid D\eu  (A_{i}) \mid 
\mbox{ : } \cup_{i=1}^{\infty} A_{i} \subset B \ , \ A_{i} \cap A_{j} = \emptyset \ \ \forall 
i \not = j \right\} \ \ .
\end{displaymath}
A function $\eu$ has bounded variation if the total variation of $D\eu$ is finite. We send the
reader to the book of Evans \& Gariepy [EG] for basic properties of such functions.

 The space $\eSBV(\Omega,R^{n})$ is defined as follows: 
$$\eSBV(\Omega,R^{n}) \ = \ \left\{ \eu \in L^{1}(\Omega,R^{n}) \mbox{ : } \mid D\eu \mid (\Omega)
 < + \infty 
\ , \ \mid D^{s}\eu \mid (\Omega \setminus \eS_{\eu}) = 0 \right\} \ .$$
The Lebesgue set of $\eu$ is the set of points where $\eu$ has approximate limit. The complementary 
set is a $\leb$ negligible set denoted by $\eS_{\eu}$. If $\eu$ is a special function with 
bounded variation then $\eS_{\eu}$ is also $\sigma$ (i.e. countably) rectifiable.  

From the Calderon \& Zygmund [CZ] decomposition theorem we obtain  the following expression of 
$D\eu$, the distributional derivative of $\eu \in \eSBV(\Omega,R^{n})$, seen as a measure: 
$$D\eu \ = \  \nabla \eu (x) \mbox{ d}x \ + \  [\eu ] \otimes \en \mbox{ d}\hen_{|_{K}} \ \ \ .$$ 

We shall use further the notation $\mu \ll \lambda$ if the measure $\mu$ is absolutely continuous 
with respect to the measure $\lambda$. 

Let us define the following Sobolev space  associated to the crack set $K$ (see [ABF]):
$$W^{1,2}_{K} \ = \ \left\{ \eu \in \eSBV(\Omega,R^{n}) \mbox{ : } \int_{\Omega} 
\mid \nabla \eu \mid^{2} 
\mbox{ d}x + \int_{K} [\eu]^{2} \mbox{ d} \hen  <  + \infty \ , \ \mid D^{s}\eu \mid \ll \hen_{|_{K}} 
\right\} \ .$$
It has been proved in [DGCL] the following equality:
\begin{equation}
W^{1,2}(\Omega \setminus K,\mathbb{R}^{n}) 
\cap L^{\infty}(\Omega,R^{n}) \ = \ W^{1,2}_{K}(\Omega,R^{n}) \cap L^{\infty}(\Omega,R^{n}) \ \ .
\label{dg}
\end{equation}
Therefore if  $\eu= \eu(\eu_{0},K)$ and $\eu_{0} \in L^{\infty}(\paromega,R^{n})$ then $\eu$  
is a special function with bounded variation.

A similar description can be made for the space of special functions with bounded
deformation $\eSBD(\Omega)$ can be found in Ambrosio, Coscia \& Dal Maso [ACDM]. 
For any function $\eu \in L^{1}(\Omega,R^{n})$ we denote by $E\eu$ the symmetric part of 
the distributional derivative of $\eu$, seen as a vector measure. We denote also by $\eJ_{\eu}$ the 
subset of $\Omega$ where $\eu$ has different approximate limits with respect to a
point-dependent direction.   The difference between $\eS_{\eu}$ and $\eJ_{\eu}$ is subtle. Let us 
quote only the fact that for a function $\eu \in \eSBV(\Omega,R^{n})$ the difference of these sets 
is $\hen$-negligible.

The definition of $\eSBD(\Omega)$ 
is the following: 
$$\eSBD(\Omega,R^{n}) \ = \ \left\{ \eu \in L^{1}(\Omega,R^{n}) \mbox{ : } \mid E\eu \mid (\Omega)
 < + \infty 
\ , \ \mid E^{s}\eu \mid (\Omega \setminus \eJ_{\eu}) = 0 \right\} \ .$$
If $\eu$ is a special function with bounded deformation then $\eJ_{\eu}$ is countably rectifiable. 
We have a decomposition theorem for $\eSBD$ functions, similar to Calderon \& Zygmund result
applied  for 
$\eSBV$ functions. The decomposition theorem is due to Belletini, Coscia \& Dal Maso [BCDM] 
and asserts that 
$$E\eu \ = \  \epsilon (\eu) (x) \mbox{ d}x \ + \  [\eu ] \odot \en \mbox{ d}\hen_{|_{\eJ_{\eu}}} \ \ \ .$$ 
Here $\odot$ means the symmetric part of tensor product and $\epsilon (\eu)$ is the approximate 
symmetric gradient, hence the approximate limit of the symmetric part of the gradient of $\eu$. 

\vspace{.5cm}

In order to give  weak versions of theorems 4.1. and 5.2. let us weaken first the space $M$ of 
pairs displacement-crack. We introduce the new set of weak pairs displacement-crack $\mathcal{ M}$: 
\begin{equation}
\mathcal{ M} \ = \  \left\{ (\eu,K) \mbox{ : } K \mbox{ is }\sigma\mbox{-rectifiable, } \eu \in
\eSBD(\Omega) \mbox{ and }  
\mid E^{s}\eu \mid (\Omega \setminus K) = 0 \right\} \ \ . 
\label{weakm}
\end{equation}
Given $(\eu,K) \in \mathcal{ M}$, the set $K$ is countably rectifiable but it is not necessarily closed; 
we impose also weaker conditions on the regularity of the displacement $\eu$. A direct 
consequence of (\ref{dg}) is that  any pair 
displacement-crack $(\eu,K)$ such that $\eu \in L^{\infty}(\Omega,R^{n})$  
belongs to the set $\mathcal{ M}$. 

Let us define the functional $\mathcal{ J}$, the weak version of the functional $J$ introduced at 
definition 4.1.: 
$$\mathcal{ J} \ :\ \mathcal{  M} \times\mathcal{ M} \rightarrow\mathbb{R}\ , \ \ \mathcal{ J} \left( 
(\eu,K),(\ev,L)\right) \ = \ 
\int_{\Omega} w(\epsilon( \ev) ) \mbox{ d}x \ + \ G \hen(L \setminus K) \ \ .$$

Before we introduce the correspondent of the function $\Psi$ from the same definition, let us 
explain what we mean by $\eu = \eu_{0}$ on the boundary of $\Omega$. We consider, for simplicity, 
that $\eu_{0}: \paromega \rightarrow\mathbb{R}^{n}$ is a continuous and therefore bounded function. Then, 
for any $\eu \in \eSBD(\Omega)$, $\eu = \eu_{0}$  if the approximate limit of $\eu$ equals 
$\eu_{0}$ in any point of $\paromega$  where the first exists, i.e.: 
$$\forall x \in \paromega, \mbox{ if  } \exists  \ev(x) \mbox{ such that } $$ 
$$lim_{\rho \rightarrow O_{+}} \frac{\int_{ B_{\rho}(x) \cap \Omega} \mid \eu(y) -  \ev(x) \mid 
\mbox{ d}y }
{\mid B_{\rho}(x) \cap \Omega \mid} = 0 \ \ \mbox{ then } \ev(x) \ = \eu_{0}(x) \ \ .$$ 
  
Let us consider a curve of imposed displacements $\lambda \mapsto \eu_{0}(\lambda) \in 
C(\paromega,R^{n})$. 
The function $\mathcal{ \Psi}$, introduced instead of $\Psi$,  is defined as follows: 
$$\mathcal{ \Psi} : \ [0,+\infty) \times \mathcal{ M} \rightarrow \left\{ 0,+ \infty \right\} \ ,$$ 
$$\mathcal{ \Psi}(\lambda,(\eu,K)) \ = \ \left\{ \begin{array}{ll}
0 & if \ \  \eu \ = \ \eu_{0} \ and \ \ \hen(K \setminus \eJ_{\eu}) = 0 \\
+\infty & otherwise \ \ . 
\end{array} \right. $$

\vspace{.5cm}

{\bf Definition 6.1.} (weak version of definition 4.1.) {\it Let us consider the space 
$\mathcal{ M}$ endowed with the topology given by the convergence: 
$$(\eu_{h},K_{h}) \ \rightarrow \ (\eu,K) \ \ if \ \ \left\{ \begin{array}{l}
\eu_{h} \ \  L^{2} \ \rightarrow \ \ \eu \ \ , \\
\hen(K_{h} \Delta K) \rightarrow 0 \ \ . 
\end{array} \right. $$
Let us consider also the function $\mathcal{ J}$, the curve of imposed displacements $t \mapsto 
\eu_{0}(t)$ with the associated function $\mathcal{ \Psi}$ and the initial data $(\eu_{0},K) \in M$ 
such that $\eu_{0} \ = \ \eu(\eu_{0}(0),K)$. 

For any $s \geq 1$ we recursively define $(\eu^{s}, K^{s}): N \rightarrow \mathcal{ M}$ like this: 

\hspace{1.cm} i) $(\eu^{s}, K^{s})(0) \ = \ (\eu_{0},K)$ ;

\hspace{1.cm} ii) for any $k \in N$  $(\eu^{s}, L^{s})(k+1) \in \mathcal{ M}$  minimizes the functional 
$$ (\ev,L) \in \mathcal{ M} \ \mapsto \ \mathcal{ J}\left( ( (\eu^{s},K^{s})(k), (\ev,L) \right) \ + \ 
\mathcal{ \Psi}((k+1)/s,
(\ev,L))$$ 
over $\mathcal{ M}$. In order to verify the constraint (\ref{cevol}), $K^{s}(k+1)$ is defined by the 
formula: 
\begin{equation}
K^{s}(k+1) \ = \ K^{s}(k) \cup \eJ_{\eu^{s}(k+1)} \ \ .
\label{gr}
\end{equation}

An energy minimizing movement associated to $\mathcal{ J}$ with the constraints (\ref{cevol}), 
$\mathcal{ \Psi}$ and 
initial data $(\eu_{0},K)$  is 
any $(\eu,K): [0,+\infty) \rightarrow \mathcal{ M}$ having the property:  
there is a diverging sequence $(s_{i})$ such that for any $t > 0$ 
$$(\eu^{s_{i}},K^{s_{i}})([s_{i}t]) \ \rightarrow (\eu,K)(t) \  \mbox{ as } i  \rightarrow \infty \
\ .$$}

\vspace{.5cm}

Let us remark that the disappearance of the set $L^{s}(k+1)$ from the crack-growth condition 
(\ref{gr}) is only apparent, because if $(\eu^{s}, L^{s})(k+1)$ minimizes the functional 
$$ (\ev,L) \in \mathcal{ M} \ \mapsto \ \mathcal{ J}\left( ( (\eu^{s},K^{s})(k), (\ev,L) \right) \ + \ 
\mathcal{ \Psi}((k+1)/s,
(\ev,L))$$ then $\mathcal{ \Psi}((k+1)/s,(\eu^{s}, L^{s})(k+1)) \ = \ 0$, hence 
$$\hen(K \setminus \eJ_{\eu}) = 0 \ \ .$$

In [ACDM] has been proven that functionals like $\mathcal{ J}$ are $L^{1}$ inferior semi-continuous 
and coercive, hence on closed subspaces ${\eV}$ of $\eSBD(\Omega)$ the functional 
$$\ev \in \eV \mapsto \mathcal{ J}\left( ( (\eu^{s},K^{s})(k), (\ev,\eJ_{\ev}) \right)$$
 has a minimizer.  Such a closed subspace of $\eSBD(\Omega)$ is the space of all weak 
displacements $\ev$ with $\ev=\eu_{0}$, where $\eu_{0}$ is a given boundary displacement. Therefore 
the following theorem is true by a trivial induction: 

\vspace{.5cm}

{\bf Theorem 4.1.}(weak version) {\it Let $\Omega \subset\mathbb{R}^{n}$ be a bounded open set with 
piecewise smooth
boundary, let $(\eu_{0},K)$ be a given admissible pair displacement-crack in $\Omega$ and let 
$$\eu_{0}: N \ \rightarrow C(\paromega,R^{n})$$ 
be a given sequence of imposed
displacements such that $\eu_{0} \ = \ \eu(\eu_{0}(0),K)$ on $\paromega \setminus K$. 

Then there exists the sequence $(\eu,K): N \ \rightarrow \ \mathcal{ M}$ such that: 

\hspace{1.cm} i) $\eu(0) \ = \ \eu_{0}$ and $K(0) \ = \ K$; 

\hspace{1.cm} ii)  for any $k \in N$ there is a countably rectifiable set $L(k+1)$ such that 
$(\eu(k+1),L(k+1)) \in \mathcal{ M}$, $\eu(k+1) \ = \ \eu_{0}(k+1)$ on $\paromega$ and 
$(\eu(k+1),L(k+1))$ is a minimizer of the functional 
$$(\ev, L) \in \mathcal{ M} \ , \ev \ = \ \eu_{0}(k+1) \mbox{ on } \paromega \ \ \mapsto 
\mathcal{ J} \left( (\eu(k),K(k)), \ev,L) \right) \ \ .$$ 
The set $K(k+1)$ is given by the formula}  $$K(k+1) \ = \ K(k) \ \cup \ \eJ_{\eu(k+1)} \ \ .$$ 

\vspace{.5cm}

The weak version of theorem 5.2. is obtained in the same way. We start by relaxing the functional 
$J_{C}$ to the functional $\mathcal{ J}_{C}$: 
$$\mathcal{ J}_{C}: \mathcal{ M} \times \mathcal{ M} \rightarrow\mathbb{R}\ ,$$
$$ \mathcal{ J}_{C}\left( (\eu,K),(\ev,L)\right) \ = \ 
\int_{\Omega} w(\epsilon( \ev) ) \mbox{ d}x \ + \ \int_{L \setminus K}f_{C}(\sig(\eu), 
\nabla \eu, \en) \mbox{ d}\hen(x) \ \ .$$
Here $\sig(\eu) \ = \ \eC \epsilon(\eu)$ and $\epsilon(\eu)$ is the approximate symmetric gradient 
of $\eu$. 

We have the following definition of an energy minimizing movement associated to $\mathcal{ J}_{C}$ with 
the usual constraints: 

\vspace{.5cm}

{\bf Definition 6.2.} (weak version of definition 5.4. adapted for $\mathcal{ J}_{C}$) 
{\it Let us consider the space 
$\mathcal{ M}$ endowed with the topology given by the convergence: 
$$(\eu_{h},K_{h}) \ \rightarrow \ (\eu,K) \ \ if \ \ \left\{ \begin{array}{l}
\eu_{h} \ \  L^{2} \ \rightarrow \ \ \eu \ \ , \\
\hen(K_{h} \Delta K) \rightarrow 0 \ \ . 
\end{array} \right. $$
Let us consider also the function $\mathcal{ J}_{C}$, the curve of imposed displacements $t \mapsto 
\eu_{0}(t)$ with the associated function $\mathcal{ \Psi}$ and the initial data $(\eu_{0},K) \in M$ 
such that $\eu_{0} \ = \ \eu(\eu_{0}(0),K)$. 

For any $s \geq 1$ we recursively define $(\eu^{s}, K^{s}): N \rightarrow \mathcal{ M}$ like this: 

\hspace{1.cm} i) $(\eu^{s}, K^{s})(0) \ = \ (\eu_{0},K)$ ;

\hspace{1.cm} ii) for any $k \in N$  $(\eu^{s}, L^{s})(k+1) \in \mathcal{ M}$  minimizes the functional 
$$ (\ev,L) \in \mathcal{ M} \ \mapsto \ \mathcal{ J}_{C}\left( ( (\eu^{s},K^{s})(k), (\ev,L) \right) \ + \ 
\mathcal{ \Psi}((k+1)/s,
(\ev,L))$$ 
over $\mathcal{ M}$. In order to verify the constraint (\ref{cevol}), $K^{s}(k+1)$ is defined by the 
formula: 
$$K^{s}(k+1) \ = \ K^{s}(k) \cup \eJ_{\eu^{s}(k+1)} \ \ .$$

An energy minimizing movement associated to $\mathcal{ J}_{C}$ with the constraints (\ref{cevol}), 
$\mathcal{ \Psi}$ and 
initial data $(\eu_{0},K)$  is 
any $(\eu,K): [0,+\infty) \rightarrow \mathcal{ M}$ having the property:  
there is a diverging sequence $(s_{i})$ such that for any $t > 0$ 
$$(\eu^{s_{i}},K^{s_{i}})([s_{i}t]) \ \rightarrow (\eu,K)(t) \  \mbox{ as } i  \rightarrow \infty \
\ .$$}

\vspace{.5cm}

The assumptions i)---iv) on $f_{C}$ from the previous section allow us to apply the main existence 
result from [ACDM]  to  the functional $\mathcal{ J}_{C}$. We have therefore the following weak version 
of the theorem 5.2.: 

\vspace{.5cm}

{\bf Theorem 5.2.}(weak version) {\it Let $\Omega \subset\mathbb{R}^{n}$ be a bounded open set with 
piecewise smooth
boundary, let $(\eu_{0},K)$ be a given admissible pair displacement-crack in $\Omega$ and let 
$$\eu_{0}: N \ \rightarrow C(\paromega,R^{n})$$ 
be a given sequence of imposed
displacements such that $\eu_{0} \ = \ \eu(\eu_{0}(0),K)$ on $\paromega \setminus K$. 

Then there exists the sequence $(\eu,K): N \ \rightarrow \ \mathcal{ M}$ such that: 

\hspace{1.cm} i) $\eu(0) \ = \ \eu_{0}$ and $K(0) \ = \ K$; 

\hspace{1.cm} ii)  for any $k \in N$ there is a countably rectifiable set $L(k+1)$ such that 
$(\eu(k+1),L(k+1)) \in \mathcal{ M}$, $\eu(k+1) \ = \ \eu_{0}(k+1)$ on $\paromega$ and 
$(\eu(k+1),L(k+1))$ is a minimizer of the functional 
$$(\ev, L) \in \mathcal{ M} \  \ , \ev \ = \ \eu_{0}(k+1) \mbox{ on } \paromega \ \ \mapsto 
\mathcal{ J}_{C} \left( (\eu(k),K(k)), \ev,L) \right) \ \ .$$ 
The set $K(k+1)$ is given by the formula}  $$K(k+1) \ = \ K(k) \ \cup \ \eJ_{\eu(k+1)} \ \ .$$ 

\vspace{.5cm}

In the anti-plane case we have to consider the space $\eSBV(\Omega,R)$ instead of 
 $\eSBD(\Omega,R^{n})$. In this case the partial regularity results of De Giorgi, Carriero \& 
Leaci [DGCL] and Ambrosio [A...] tell us that the classical Mumford-Shah functional has minimizers 
in the set of pairs displacement-crack $M$. Both theorems are therefore true in the strong form, 
in the anti-plane case.

\subsection{Theorem 4.2.} 
\indent

For any $D > 0$ let $u_{0}(D)$ be the following boundary displacement: 
$$u_{0}(D)(x) \ = \ \left\{ \begin{array}{ll}
0 & on \ \ \gau^{1} \\
D & on \ \ \gau^{2} \ \ . 
\end{array} \right. $$
For any crack set $K$ we denote by $u_{K}$ the displacement $u_{K} \ = \ u(u_{0}(1),K)$, i.e. the 
solution (or one of the solutions) of the problem: 
$$ \left\{ \begin{array}{ll}
         div \nabla u = 0 & in \ \  \Omega \setminus K \\
         \nabla u \en = 0 & on \ \  \gaf \cup K \\
         u = 0 & on \ \  \gau^{1} \\
         u = 1 & on \ \  \gau^{2} \ \ . 
         \end{array}  \right. $$
More general, we shall use the notation $u_{K}(D) \ = \ u(u_{0}(D),K)$. It is obvious that 
$u_{K}(D) \ = \ D \ u_{K}$. 
In the body of Theorem 4.2. we have introduced the Mumford-Shah functional: 
$$I(v, K) \ = \ \frac{1}{2} \ \int_{\Omega} \mid \nabla v \mid \mbox{ d}x \ + \ G \mathcal{ H}^{1}(K 
) \ \ .$$
The displacement $u_{K}$ has the minimum property 
$$I(u_{K},K) \ \leq \ \ I(v,K) \ \ ,  \ \ \forall (v,K) \in M  \ , \ v \ = \ u_{0}(1) \mbox{ on } 
\gau^{1} \cup \gau^{2} \ \ .$$
That is why it is reasonable to redefine the functional $I$ as a functional depending only on the
crack set $K$: 
$$\tilde{I}(K) \ = \ I(u_{K},K) \ \ .$$ 
With this notation we have, for any $D > 0$, the inequality: 
$$D^{2} \  \ \frac{1}{2} \ \int_{\Omega} \mid \nabla u_{K} \mid \mbox{ d}x \ + \ G \mathcal{ H}^{1}(K 
) \ \leq \ I(v,K) \ \ ,  \ \ \forall (v,K) \in M  \ , \ v \ = \ u_{0}(D) 
\mbox{ on } \gau^{1} \cup \gau^{2} \ \ .$$
We make the notation: 
$$\tilde{I}(K,D) \ = \  D^{2} \  \ \frac{1}{2} \ \int_{\Omega} \mid \nabla u_{K} \mid \mbox{ d}x \ + \ G \mathcal{ H}^{1}(K 
) \ \ .$$
 
We shall need further the function $\ueo$, which is defined modulo an additive
constant by the relations: 
$$\frac{\partial \ueo}{\partial x_{1}} = \frac{\partial u_{\emptyset}}{\partial x_{2}} \ \ ,$$
$$\frac{\partial \ueo}{\partial x_{2}} = - \frac{\partial 
u_{\emptyset}}{\partial x_{1}} \ \ .$$
The level sets of $\ueo$ form a congruence of curves in $\overline{\Omega}$.   The part of the boundary 
$\gaf$ belongs to this congruence. We define the following system of open neighbourhoods named 
$V(\ueo)$,  with 
the aid of this congruence: 
$$\forall A \in V(\ueo) \mbox{ \ \ \ }  \partial A \  \setminus \gau \mbox{ 
it is locally a level set of }
 \ueo \ \ .$$
For any $A \in V(\ueo)$ we denote by $\partial_{u} A$ the part of the boundary of $A$ belonging to 
$\gau^{1}$ or $\gau^{2}$, i.e. 
$$\partial_{u} A \ = \ \partial A \ \cap  \ \left( \gau^{1} \cup \gau^{2} \right) \ \ .$$
The remaining part of $\partial A$ is denoted by $\partial_{f} A$. 

Let $K$ be a rectifiable curve and $\Omega' \in V(\ueo)$ such that 
$$K \setminus \gaf  \ \subset \ \Omega'  \cup \partial_{\eu} \Omega' \ \ .$$ 
For the couple $(K,\Omega')$ we introduce the following stress field: 
$$ \sigma = \left\{ \begin{array}{ll}
                     \nabla u_{\emptyset} & in \ \ \Omega \setminus \Omega' \\
                       0      & in \ \ \Omega' \ \ . 
                    \end{array} \right. $$
This stress field is statically admissible with respect to the body with reference configuration 
$\Omega \setminus K$ and boundary displacement $u_{0}(1)$. Therefore we have the following 
inequality: 
$$\frac{1}{2} \int_{\Omega} \mid \nabla u_{K} \mid^{2} \mbox{ d}x \geq \int_{\gau} (\sigma \en) 
\cdot u_{0}(1) \mbox{ d}\mathcal{ H}^{1} - \frac{1}{2} \int_{\Omega} \mid \sigma \mid^{2} \mbox{ d}x \ \ 
.$$
The latter inequality can be put in terms of Mumford-Shah functional $I$ like this: 
$$\tilde{I}(\emptyset) - \tilde{I}(K) \leq \  \frac{1}{2}\int_{\Omega'} \mid \nabla u_{\emptyset} 
\mid^{2} \mbox{ d}x \  -  \  G \mathcal{ H}^{1}(K 
\setminus \gaf) \ \ .$$
The reason for which we have put $\mathcal{ H}^{1}(K 
\setminus \gaf)$ instead of $\mathcal{ H}^{1}(K)$ is that $u_{K} \ = \ u_{K \setminus \gaf}$ but 
$\tilde{I}(K) \geq \tilde{I}(K \setminus \gaf)$. As a consequence, sets $K$ with a part on $\gaf$ 
are disqualified to be minimizers of $\tilde{I}$. 

Let us denote by $\tau$ the tangent vector field in direct sense to $\paromega'$. After few
calculations we obtain from the previous inequality the estimation: 
\begin{equation}
\tilde{I}(\emptyset) - \tilde{I}(K) \leq \frac{1}{2} \int_{\paromega'} u_{\emptyset} 
\cdot (
\nabla \ueo \tau ) - \ G \mathcal{ H}^{1}(K 
\setminus \gaf) \ \ .
\label{fest}
\end{equation}
We deduce that for any $D>0$ we have:
\begin{equation}
\tilde{I}(\emptyset,D) - \tilde{I}(K,D) \leq \frac{D^{2}}{2} \int_{\paromega'} u_{\emptyset} 
\cdot (
\nabla \ueo \tau ) - \ G \mathcal{ H}^{1}(K 
\setminus \gaf) \ \ .
\label{fest1}
\end{equation}

Let us return to the congruence of curves defined by $\ueo$ and consider the projection function 
on $\gau^{2}$ with respect to the congruence. For any set $B \subset \gau^{2}$ we denote by 
$V \ueo (B)$ the variation of $\ueo$ on $B$. We see that: 
\begin{equation}
\frac{D^{2}}{2} \int_{\partial_{u} \Omega^{'}}  (\nabla u_{\emptyset} \en) \cdot 
u_{\emptyset} \mbox{ d}
\mathcal{ H}^{1} \geq \frac{1}{2} D^{2} \ V \ueo (P(K)) \ \ ,
\label{l1}
\end{equation}
because of the equality: 
$$\inf \left\{ \frac{D^{2}}{2} \int_{\partial_{u} \Omega^{'}}  (\nabla u_{\emptyset} \en) 
\cdot u_{\emptyset} 
\mbox{ d}\mathcal{ H}^{1} \mbox{ : } \Omega^{'} \in V(\ueo), K \setminus \gaf \subset \Omega^{'}
\right\} = 
\frac{1}{2} D^{2} \ V \ueo (P(K)) \ \ .$$
From the inequalities (\ref{fest1}) and (\ref{l1}) we obtain the improved estimation:
\begin{equation}
 \tilde{I}(\emptyset,D) - \tilde{I}(K,D) \ \leq \  \frac{1}{2} D^{2} \ V \ueo (P(K)) \ - \  
G \mathcal{ H}^{1}(K \setminus \gaf) \ \ .
\label{l2}
\end{equation} 
Let us remark that $$V\ueo(K) \geq V\ueo(P(K))$$ therefore we have: 
$$ \tilde{I}(\emptyset,D) - \tilde{I}(K,D) \ \leq \  \frac{1}{2} D^{2} \ V \ueo (K) \ - \  
G \mathcal{ H}^{1}(K \setminus \gaf) \ \ .$$
Due to the assumption (recall the notation $u(D) = u(u_{0}(D),\emptyset)$) 
\begin{equation}
C \geq \mid \nabla u(1)\en \mid^{2} \geq c > 0
\label{bound}
\end{equation}
we have $ V \ueo (K) \leq C \mathcal{ H}^{1}(K)$ hence (if we suppose that $K \cap \gaf$ is 
$\mathcal{ H}^{1}$ negligible): 
\begin{equation}
\tilde{I}(\emptyset,D) - \tilde{I}(K,D) \ \leq \ \left[  \frac{1}{2} D^{2} C  \ - G \right] 
\mathcal{ H}^{1}(K) \ \ .
\label{e1}
\end{equation}
Therefore, if 
$$\frac{1}{2} D^{2} C  \ < \ G \ \ ,$$
from (\ref{e1}) we see that $\emptyset$ minimizes  $\tilde{I}(\cdot,D)$, which proves the point i) 
of the theorem. 
Let us go back to (\ref{l2}) and introduce $\Gamma(K)$ as the curve  with the properties: 

\hspace{2.cm} p1) for any $\Omega' \in V(\ueo)$, if $K \setminus \gaf \subset \Omega'$ then 
$\Gamma(K) \subset \Omega'$, 

\hspace{2.cm} p2) $\Gamma(K)$ is a length minimizer in the class of curves that fulfills p1).

We remark that $\Gamma(K)$ might not be unique, but it always exists.

It is straightforward that $\Gamma(K) = \Gamma(P(K))$ and $V \ueo (P(K)) = V \ueo (\Gamma(K))$. 
We have then: 
\begin{equation}
 \tilde{I}(\emptyset,D) - \tilde{I}(K,D) \ \leq \  \frac{1}{2} D^{2} \ V \ueo (\Gamma(K)) \ - \  
G \mathcal{ H}^{1}(\Gamma(K)) \ \ .
\label{e2}
\end{equation} 
From the   assumption (\ref{bound})  we see that
\begin{equation}  
\frac{1}{2} D^{2} \ V \ueo (P(K)) \ - \  
G \mathcal{ H}^{1}\Gamma(P(K)) \ \geq \ \left[ \frac{D^{2}}{2} c \ - G \right] \mathcal{ H}^{1}(\Gamma(K)) \ 
\ .
\label{e3}
\end{equation} 
Therefore, if $$\frac{D^{2}}{2} c \ >  \ G $$ then the right member of (\ref{e3}) is positive and 
it attains the maximum when $\mathcal{ H}^{1}(\Gamma(K))$ is maximal. This happens when $\Gamma(K)$ 
separates $\gau^{1}$ from $\gau^{2}$. In this case is easy to see that we have equality in 
the relation (\ref{e2}), which proves the point ii) of the theorem.  

The proof of iii) it is now straightforward. If $C=c$ then 
$$\frac{1}{2} D^{2} C \ = \   \frac{D^{2}}{2} c  \ \ .$$

\section{Conclusions and perspectives}
\indent

The first model contains only a constant connected to fracture, namely the constant of 
Griffith $G$. The main qualities of this model are: 

\hspace{2.cm} i) crack appearance is allowed, together with 
crack propagation, 

\hspace{2.cm} ii) there is no restriction concerning the pattern of the crack during its 
evolution.

We have seen that in the first model the critical stress which lead to fracture (or crack 
appearance) is not a constant of material. 

The second model contains two constants of material connected to fracture: $G$ and a constant with 
the dimension of a stress named $\Sigma$. In this model the critical stress which lead to crack 
appearance is a constant of material, related to $\Sigma$. This model has the same qualities as the 
first. 

These two models are fully macroscopical, in the sense that no fracture mechanism based on 
micro-cracks or other micro-defects  was supposed.

The main open theoretical problem is the general existence of an energy minimizing movement according 
to our definitions. Below is described an existence result based on a sound physical 
assumption (\ref{power}). Nevertheless, we do not know if  (\ref{power}) can be proved from the 
basic assumptions of the model.    

\vspace{.5cm}

{\bf Theorem 7.1.} {\it 
Let us consider for a given $s$ an incremental solution $k \mapsto 
(\eu^{s}(k),K^{s}(k)) \in M$, according to definition 4.1.. For any $k \in N$ we introduce the 
displacement 
$$\ev^{s}(k+1) \ = \ \eu( (k+1)/s,K^{s}(k)) \ \ .$$
Let us suppose that the power communicated by the rest of the universe to the body is bounded at any 
moment $t$.  
The incremental form of this assumption  consists in the existence of a constant $P$ such that  
for any $k$ and $s$ we have 
\begin{equation}
\langle \eT(K^{s}(k))\frac{1}{2}\left( \eu_{0}((k+1)/s) + \eu_{0}(k/s) \right) , 
\eu_{0}((k+1)/s) - \eu_{0}(k/s) \rangle  \ \leq \ P/s \ \ .
\label{power}
\end{equation}

Then for any $t > 0$ there 
exist diverging sequences $(s_{i})_{i}$ and $(k_{i})_{i}$ such that $k_{i}/s_{i}$ converges to 
$t$ and $(\eu^{s_{i}}, K^{s_{i}})(k_{i})$ converges to an element of $M$ $(\eu,K)(t)$. }

\vspace{.5cm}
{\bf Proof: } 
From the minimality assumption on the incremental solution we have for any $k \in N$  
the inequality: 
$$J((\eu^{s}(k),K^{s}(k)),(\ev^{s}(k+1),K^{s}(k))) \ \geq \ J((\eu^{s}(k),K^{s}(k)),
(\eu^{s}(k+1),K^{s}(k+1))) \ \ .$$
This inequality means that: 
$$\int_{\Omega} w(\nabla \ev^{s}(k+1)) \mbox{ d}x \  \geq  \ 
\int_{\Omega} w(\nabla \eu^{s}(k+1)) \mbox{ d}x \ + \ G \hen(K^{s}(k+1) \setminus K^{s}(k)) \ \ .$$
The crack growth condition $K^{s}(k) \subset K^{s}(k+1)$ implies that the latter relation  
can be put in the following form: 
\begin{equation}
\left( \int_{\Omega} w(\nabla \ev^{s}(k+1)) \mbox{ d}x \ - \ 
\int_{\Omega} w(\nabla \eu^{s}(k)) \mbox{ d}x  \right) \ + \ 
\label{inc}
\end{equation} 
$$ + \ \int_{\Omega} w(\nabla \eu^{s}(k)) \mbox{ d}x   \ +  \ G \hen(K^{s}(k) \ \geq 
 \ \int_{\Omega} w(\nabla \eu^{s}(k+1)) \mbox{ d}x \ + \ G \hen(K^{s}(k+1)) \ \ .$$ 
This is the incremental form of the Griffith criterion of crack  propagation (\ref{grif}). Indeed, 
we have the chain of  equalities: 
$$\int_{\Omega} w(\nabla \ev^{s}(k+1)) \mbox{ d}x \ - \ 
\int_{\Omega} w(\nabla \eu^{s}(k)) \mbox{ d}x \ = $$ 
$$ = \ \frac{1}{2} \langle \eT(K^{s}(k)) \eu_{0}((k+1)/s), 
 \eu_{0}((k+1)/s) \rangle  \ - \ \frac{1}{2} \langle \eT(K^{s}(k)) \eu_{0}(k/s),\eu_{0}(k/s) \rangle \ = $$ 
$$ = \ \langle \eT(K^{s}(k))\frac{1}{2}\left( \eu_{0}((k+1)/s) + \eu_{0}(k/s) \right) , 
\eu_{0}((k+1)/s) - \eu_{0}(k/s) \rangle \ \ .$$ 
$\ev^{s}(k+1)$ represents the displacement of the body with the boundary displacement 
$\eu_{0}(k/s + 1/s)$ in the presence of the crack $K^{s}(k)$. $\eu^{s}(k)$ represents 
the displacement of the body with the boundary displacement 
$\eu_{0}(k/s)$ in the presence of the same crack $K^{s}(k)$. According to (\ref{extp}), the
quantity 
$$ \left( \int_{\Omega} w(\nabla \ev^{s}(k+1)) \mbox{ d}x \ - \ 
\int_{\Omega} w(\nabla \eu^{s}(k)) \mbox{ d}x \right) / \left( \frac{1}{s} \right)$$ 
is the discretized expression of the power communicated by the rest of the universe to the body 
at the moment $k/s$, when a time discretization with step $1/s$ is considered.

We deduce from the inequality (\ref{inc}) that 
$$ P/s \ + \ \int_{\Omega} w(\nabla \eu^{s}(k)) \mbox{ d}x   \ +  \ G \hen(K^{s}(k) \ \geq 
 \ \int_{\Omega} w(\nabla \eu^{s}(k+1)) \mbox{ d}x \ + \ G \hen(K^{s}(k+1)) \ \ .$$ 
We have therefore: 
$$P k/s \ \geq \  \int_{\Omega} w(\nabla \eu^{s}(k+1)) \mbox{ d}x \ + \ G \hen(K^{s}(k+1)) \ \ .$$
From the compactness theorem for $\eSBD$ space and the latter inequality we deduce that for
any $t > 0$ there 
exist diverging sequences $(s_{i})_{i}$ and $(k_{i})_{i}$ such that $k_{i}/s_{i}$ converges to 
$t$ and $(\eu^{s_{i}}, K^{s_{i}})(k_{i})$ converges to an element of $M$ $(\eu,K)(t)$. 
\hfill $\Box$

\vspace{.5cm}

 In the paper [AB] Ambrosio \& Braides introduce a generalized minimizing 
movement based model for the propagation of a crack in the presence of viscous forces in the body. 
They give as initial datum at $t=0$ the anti-plane displacement $u_{0} \in \eSBV(\Omega,R) \cap 
L^{\infty}(\Omega,R)$. 
For a given $s$ they recursively define a sequence $(u^{s}_{k})_{k}$ in $\eSBV(\Omega,R)$ and 
an increasing sequence of closed rectifiable sets $(K^{s}_{k})_{k}$ as follows: $u_{0}^{s} =
 u_{0}$,  $K^{s}_{0}=\emptyset$ and $u^{s}_{k+1}=w$, $K^{s}_{k+1}= \overline{\eS_{w}} \cup 
K^{s}_{k}$, where $w$ is a minimizer of the functional 
\begin{equation}
v \mapsto  \int_{\Omega}  \mid \nabla v \mid^{2} \mbox{ d}x + \hen(\eS_{v} \setminus
K^{s}_{k}) + s \int_{\Omega} \mid v - u^{s}_{k} \mid^{2} \mbox{ d}x 
\label{ab}
\end{equation} 
over the set of all $v$ such that: 
$$ v \in \eSBV(\Omega,R) \ , \  \| v \|_{\infty} \leq \| u_{0} \|_{\infty}  \ \ .$$ 
The generalized minimizing movements obtained as limits of such incremental solutions, when 
$s$ diverges, correspond to the following situation: a body evolves from the initial state $u_{0}$, 
with the initial crack $\eS_{u_{0}}$, under a constant imposed boundary displacement. The equation
of evolution for the displacement is: 
$$div \ \nabla u(t) \ + \ \dot{u} (t) \ = \ 0 \ \ .$$

The authors obtain an existence result for the generalized minimizing movement introduced by them. 
After the introduction of the piecewise constant function: 
$$u^{s}(t) \ = \ u^{s}_{[st]} \ \ ,$$ 
they find the following estimation: 
\begin{equation}
\| u^{s}(t') - u^{s}(t) \|_{L^{2}} \ \leq \ M \sqrt{t' - t + \frac{1}{s}} \ \ if \  t' \geq t \ \ . 
\label{est}
\end{equation}
Therefore there exists a diverging sequence $(s_{i})_{i}$ such that $u^{s_{i}}$ converges to $u$ 
uniformly in $L^{\infty}([0,T], L^{2}(\Omega,R))$, for all $T > 0$. A consequence of this result 
is that the crack appearance is forbidden in this model. 

This result is obtained under the assumption of constant imposed boundary displacement, 
 equal to the  
trace on the boundary of the initial datum $u_{0}$. 

It is natural to introduce the Lam\'e constant $\mu$ and the viscosity $\lambda$ in the expression 
of the functional (\ref{ab}) and modify it like this: 
$$v \mapsto  \int_{\Omega} \mu  \mid \nabla v \mid^{2} \mbox{ d}x + \hen(\eS_{v} \setminus
K^{s}_{k}) + \lambda s \int_{\Omega} \mid v - u^{s}_{k} \mid^{2} \mbox{ d}x \ \ .$$ 
We obtain the more physical case of an anti-plane displacement satisfying at any moment 
$t$ the equation: 
$$div \ \mu \nabla u(t) \ + \ \lambda \dot{u} (t) \ = \ 0 \ \ .$$
The estimation (\ref{est}) becomes 
$$\| u^{s}(t') - u^{s}(t) \|_{L^{2}} \ \leq \ M \sqrt{t' - t + \frac{1}{\lambda s}} \ \ if \  t' 
\geq t \ \ .$$
We expect to obtain our first model, in the case of anti-plane displacements,  when the viscosity 
$\lambda$ converges to $0$. It is easy to see that if $\lambda$ converges to $0$ then the uniform 
estimation from above is lost, hence there is no contradiction between the fact that in our model 
crack appearance is allowed and the fact that in the model of Ambrosio \& Braides crack
appearance is forbidden.     

As a conclusion, an open direction of research consists in the use of more general minimizing 
movements 
in order to study the propagation of a crack in the presence of viscous effects (as is the paper 
[AB]) or in the case 
of an elasto-plastic body.

The models presented in the paper are of applicative interest. In order to use them we have 
to know how to minimize a Mumford-Shah functional. This can be done by approximating, in the sense 
of variational convergence, the original functional by a less strange one. The idea is to 
replace the pair displacement-crack  $(\eu,K)$ with the pair $(\eu,f)$, where $f$ is a 
smoothed version of the characteristic function of the crack set $K$, taking values in the 
interval $[0,1]$. The original functional may be replaced by an Ambrosio-Tortorelli approximation, 
introduced in [AT1], [AT2]. This opens the path to future interesting numerical results.

\newpage

\centerline{{\large References}}

\vspace{.5cm}

\begin{itemize}

\item
[A1] L. Ambrosio, Variational problems in SBV and image segmentation, Acta Appl.  
Mathematic\ae  17, 1989,1-40

\item
[A2] L. Ambrosio, Existence Theory for a New Class of Variational Problems, 
Arch. Rational Mech. Anal., vol. 111, 1990, 291-322

\item
[A3] L Ambrosio, The space $SBV(\Omega)$ and free discontinuity problems, in 
Variational and Free Boundary Problems, editori A. Friedman, J. Spruck, IMA Vol. 
in Math. and Its Appl., vol. 53, Springer-Verlag, 1994, 1-24

\item
[AB] L. Ambrosio, A. Braides, Energies in SBV and Variational Models in Fracture Mechanics, 
Proceedings of the EurHomogenization congress, Nizza,   
Gakuto Int. Series, Math. Sci. and Appl.,  9, 1--22,
1997.

\item
[ABF] L. Ambrosio, G. Buttazzo, I. Fonseca, Lower semicontinuity problems in Sobolev spaces 
with respect to a measure, J. Math. Pures Appl. 75, 1996, 211-224

\item
[ACDM] L. Ambrosio, A. Coscia, G. Dal Maso, Fine Properties of Functions with Bounded Deformation, 
Preprint SISSA 8/96/M, 1996

\item
[All] W.K. Allard, On the first variation of a varifold, Ann. of Math., vol. 95, no. 3, 1972, 
417-491  

\item
[AT1] L. Ambrosio, V.M. Tortorelli, Approximation of Functionals Depending on Jumps by Elliptic Functionals via 
$\Gamma$-convergence, Comm. Pure Appl. Math., vol. 43, 1990, 999-1036

\item
[AT2]  L. Ambrosio, V.M. Tortorelli, On the approximation of free discontinuity problems, 
Boll. U.M.I., 6-B, 1992, 105-123

\item
[BCDM] G. Bellettini, A. Coscia, G. Dal Maso, Compactness and lower semicontinuity properties 
in $\eSBD(\Omega)$, Preprint S.I.S.S.A. 86/96/M, 1996

\item
[Ba] J. Ball, Some recent developments in nonlinear elasticity and its applications to 
material sciences, to appear in Proc. EPSRC Spring School 1995, Cambridge Univ. 
Press

\item
[Bu1] M. Buliga, Variational Formulations in Brittle Fracture Mechanics,  PhD Thesis, 
Institute of Mathematics of the Romanian Academy, 1997

\item
[Bu2] M. Buliga, Modelisation de la d\'ecohesion d'interface fibres-matrice dans les materiaux composites, 
memoire de D.E.A., Ecole Polytechnique, 1995

\item
[Bu3] M. Buliga, Energy concentration and brittle crack propagation, J. of Elasticity, 52, 3, 201-238, 1999

\item
[CZ] A.P. Calderon, A. Zygmund, On the differentiability of functions which are of bounded
variation in Tonelli's sense, Rev. Un. Mat. Argentina, 20, 1960, 102-121

\item
[DGA] E. De Giorgi, L. Ambrosio, Un nuovo funzionale del calcolo delle variazioni, Atti Accad. Naz. 
Lincei Rend. Cl. Sci. Fis. Mat. Natur., 82, 1988, 199-210

\item
[DGCL] E. De Giorgi, G. Carriero, A.  Leaci, Existence theory for a minimum problem with free 
discontinuity set, Arch. Rational Mech. Anal., vol. 108, 1989, 195-218 

\item
[DP] G. Del Piero, Recent developments in the mechanics of materials which do not support tension, 
in Free Boundary Problems: Theory and 
Applications, vol I, Eds. Hoffmann K. H., Sprekels J., Pitman res. notes in math. series, 
Longman Scientific \& Technical, 1990

\item
[EbM] D. G. Ebin, J. Marsden, Groups of diffeomorphisms and the motion of an incompressible fluid, Ann. Math., vol. 92, 
no.1 , 1970, 102-163

\item
[Es] J. D. Eshelby, Energy relations and the energy-momentum tensor in continuum mechanics, 
Inelastic Behavior of Solids, ed. M.F. Kanninen et al., New York: McGraw-Hill, 1970, 77-115

\item
[EG] L. C. Evans, R. F. Gariepy, Measure Theory and Fine Properties of Functions, CRC Press, 1992

\item
[FMa] G. Francfort, J.-J. Marigo, Stable damage evolution in a brittle continuous medium, 
Eur. J. Mech., A/Solids, 12, no. 2, 1993, 149-189

\item
[G] A.A. Griffith, The phenomenon of rupture and flow in solids, Phil. Trans. Royal Soc. London, A 221, 
1920, 163-198

\item
[Gu1] M. E. Gurtin, On the Energy Release Rate in Quasistatic Elastic Crack Propagation, 
J. of Elasticity, vol 9, no. 2, 1979, 187-195

\item
[Gu2] M. E. Gurtin, Thermodynamics and the Griffith Criterion for Brittle Fracture, 
Int. J. Solids Structures, vol. 15, 1979, 553-560

\item
[I] G.R. Irwin, Structural Mechanics, Pergamon Press, London, England, 1960

\item
[MS] D. Mumford, J. Shah, Optimal approximation by piecewise smooth functions and associated 
variational problems, Comm. on Pure and Appl. Math., vol. XLII, 
no. 5, 1989, 577-685

\item
[Oht1] K. Ohtsuka, Generalized J-integral and Its Applications I. Basic Theory, Japan J. Appl. Math., 
2, 1985, 21-52

\item
[Oht2] K. Ohtsuka, Generalized J-integral and three-dimensional fracture mechanics I, 
Hiroshima Math. J.,11, 1981, 329-350

\item
[Oht3] K. Ohtsuka, Generalized J-integral and three-dimensional fracture mechanics II, 
Hiroshima Math. J., 16, 1986, 327-352

\item
[Oht4] K Ohtsuka, Generalized J-integral and its applications, RIMS Kokyuroku A62, 149-165

\item
[R] J.R. Rice, Mathematical analysis in the mechanics of fracture, in Fracture: an Advanced Treatise, vol. 2, ed. 
H. Liebowitz, Academic Press, 1969, 191-311

\item
[StLe] H. Stumpf, K. Ch. Le, Variational principles of nonlinear fracture mechanics, 
Acta Mechanica 83, 1990, 25-37

\end{itemize}

\end{document}